\documentclass[11pt]{article}
\usepackage[english]{babel}
\usepackage{latexsym}
\usepackage{amsfonts,amsmath}
\usepackage{color}
\usepackage{comment}
\newsymbol\blacktriangleright 1349
\newsymbol\blacktriangleleft 134A
\newsymbol\vartriangle 134D

\newtheorem{thm}{Theorem}
\newtheorem{corol}[thm]{Corollary}
\newtheorem{lemma}[thm]{Lemma}
\newtheorem{fact}[thm]{Fact}

\newtheorem{exa}[thm]{Example}
\newtheorem{remark}[thm]{Remark}
\newenvironment{example}{\begin{exa} \rm }{\hfill$\Box$ \end{exa}}

\newenvironment{proof}{\begin{trivlist}\item[] \mbox{\it Proof. }}
{\hfill$\Box$ \end{trivlist}}

\definecolor{lgrey}{gray}{.55}%
\definecolor{hgrey}{gray}{.3}%

\def\dv{\mathbb}
\def\caP{{\cal P}(N)}
\def\enuN{\Upsilon (N)} 
\def\chain{{\cal C}}
\def\cone{\mbox{\rm con}\,}
\def\conv{\mbox{\rm conv}\,}
\def\Cov{\mbox{\rm Cover}\,} 
\def\cover{\prec\!\!\cdot\,\,}
\def\cl{\mbox{\rm cl}\,}  
\def\dist{\mbox{\rm dist}\,} 
\def\diag{\Delta} 
\def\tr{\mbox{\rm tr}} 
\def\sym{\mathfrak{S}} 
\def\permA{\varepsilon} 
\def\permB{\eta} 
\def\permC{\sigma} 
\def\permD{\pi} 
\def\permE{\zeta} 
\def\permF{\xi} 
\def\permG{\varsigma} 
\def\covA{\prec_{\varepsilon}\!\!\cdot\,\,} 
\def\perN{\Pi(N)} 
\def\preced{\circ} 
\def\braid{\bullet} 
\def\topol{\textcolor{lgrey}{\bullet}} 
\def\lGalo{\triangleleft} 
\def\rGalo{\triangleright} 
\def\lGprec{\lhd} 
\def\rGprec{\rhd} 
\def\lGbraid{\blacktriangleleft} 
\def\rGbraid{\blacktriangleright} 
\def\lGtopol{\textcolor{lgrey}{\blacktriangleleft}} 
\def\rGtopol{\textcolor{lgrey}{\blacktriangleright}} 
\def\rankv{\rho} 
\def\varrankv{\varrho} 
\def\trans{\tau} 
\def\Inv{\mbox{\rm Inv}\,} 
\def\latt{{\cal L}}  
\def\nor{{\mathsf N}} 
\def\calF{{\cal F}}
\def\calX{{\cal X}}
\def\calY{{\cal Y}}
\def\calD{{\cal D}}

\begin{document}
\addtolength{\baselineskip}{1mm}
\parskip 0mm

\title{Graphical view on linear extensions of finite posets}
\author{{\Large Milan Studen\'{y} ~and~ V\'{a}clav Kratochv\'{\i}l}\\[0.3ex]
Czech Academy of Sciences,\\
Institute of Information Theory and Automation}
\date{\today} 
\maketitle

\begin{abstract}
One of the possible cryptomorphic definitions of a partially ordered set (= a poset) $P$
on a non-empty finite ground set $N$ is in terms of the set ${\cal L}(P)$ of all its linear extensions,
that is, in terms of the set of total orders on~$N$ consistent with $P$. Any total order on $N$ can be interpreted
as a node of a particular graph, called the {\em permutohedral graph\/} (over $N$), because it is indeed
the graph of a certain polytope $\perN$ in ${\dv R}^{N}$, known as the permutohedron.

It is shown in the paper that a non-empty set of total orders on $N$ equals to ${\cal L}(P)$
for some poset $P$ on~$N$ if and only if it is a geodetically convex set in the permutohedral graph.
This result means that a purely graphical concept of geodetic convexity in this graph is a cryptomorphic definition of a finite poset. In particular, the lattice of geodetically convex sets in this graph is graded and its height function is described in graphical terms. A counter-example, however, shows that the height function does not correspond to the usual graphical diameter, relating this matter to a combinatorial concept of the
{\em dimension of a poset\/}.

Two alternative cryptomorphic views on a poset $P$ on $N$ are also 
discussed. The geometric counterpart is its full-dimensional {\em braid cone\/} in ${\dv R}^{N}$,
while a combinatorial alternative is a {\em topology\/} on $N$ distinguishing points,
often referred as a (finite) distributive lattice.
\end{abstract}


\section{Introduction}\label{sec.introduction}
The motivation of the first author for considering different interpretations of {\em partially ordered sets\/} \mbox{(= posets)}
on a non-empty finite ground set $N$ originates from his recent findings in \cite{Stu24}.
The main message there is that non-empty faces of the cone of supermodular functions on the power set
$\caP$ of $N$ have combinatorial interpretations and the most promising approach seems to be to describe
these faces by certain collections of {\em posets on $N$}, introduced already in \cite[\S\,3.3]{PRW08}
and called the {\em fans of posets}. These mathematical objects, however, can also alternatively be
interpreted as collections of particular cones in ${\dv R}^{N}$, sometimes named {\em braid cones}.
Yet another version of these objects are particular collections of {\em topologies on $N$},
which seem to correspond to what was discussed in \cite[\S\,3.3(d)]{Fuj91}.
Further equivalent mathematical objects are {\em rank tests\/} discussed in \cite{MPSSW09},
which are certain partitions of the set of all total orders on $N$. The rank tests correspond
to some statistical hypotheses tests studied in context of algebraic statistics.

In mathematics, two objects are considered to be cryptomorphic if they are equivalent, but not obviously equivalent.
To put the cryptic equivalence of above mathematical objects on a solid base, one first needs to establish formally
the equivalence between the components of those objects. In the given context, this means clarifying what are precisely
the respective equivalent views on posets on a fixed finite set $N$, which is the topic of the present paper.
In our opinion, the most elegant way to introduce the alternative views on finite posets is to use the tools
and results from lattice theory on (order-reversing) {\em Galois connections\/} \cite[\S\,V.8]{Bir95},
as done below.
\smallskip

This paper is predominantly devoted to the characterization of a poset $P$ on $N$ in terms of
the collection ${\cal L}(P)$ of all its {\em linear extensions}. The correspondence $P\leftrightarrow {\cal L}(P)$
has earlier been treated in the literature, both in mathematics and in computer science.
Pruesse and Ruskey \cite{PR94} proposed an algorithm for generating (all) linear extensions of a given
poset~$P$ and discussed the computational complexity of the tasks of either generating or counting linear extensions.
In their paper, they, like many previous or later authors, described total orders on $N$ by ordered lists of
(all) elements of $N$, which we below name {\em enumerations\/} of $N$. They also utilized
a natural proximity relation among these when two enumerations are considered to be neighbors
if one is obtained from the other by the so-called {\em adjacent transposition}.
These neighborhood relations correspond to the edges in our {\em permutohedral graph}, discussed below.

Linear extensions of a finite poset $P$ were also a topic of interest in several more recent mathematical publications.
The restriction of the permutohedral graph to ${\cal L}(P)$, called the {\em linear extension graph\/}
(of $P$), was studied from the point of view of graph theory in \cite{Naa00,Mas09}. A recent
paper \cite{CP23} offers a broad survey of inequalities for the number of linear extensions of
a poset $P$, based on tools of enumerative combinatorics. One of the classic approaches in this field
is to assign a certain $n$-dimensional {\em order polytope\/} to a poset $P$ on an $n$-element set $N$
such that the geometric structure of this polytope reflects the combinatorial structure of $P$ \cite{Sta86}.
For example, the volume of the order polytope is the normalized number of linear extensions, being
the number of linear extensions of $P$ divided by $n!$.  Stanley in his book \cite[p.\ 218-221]{Sta97} also studied
the {\em order polynomial\/} of $P$, which also reflects the combinatorial structure of $P$. In fact, the order polynomial is related to the so-called Ehrhart polynomial of the order polytope \cite{Sta86} and the leading coefficient of both these polynomials is the volume of the order polytope of $P$ (further links to enumerative combinatorics are discussed in the remarks from Section~\ref{sec.other-views}).
\smallskip

An important related concept is that of the {\em dimension of a poset $P$} on $N$, introduced by
Dushnik and Miller \cite{DM41} already in the 1940s. Briefly formulated, it is the smallest
number of total orders on $N$ whose intersection is $P$. A classic result \cite{Hir51} in the theory of
poset dimension says that the dimension of a poset $P$ on an $n$-element set $N$, $n\geq 4$, is at most
$\lfloor \frac{n}{2}\rfloor$; nonetheless, this upper bound is tight owing to a classic construction
of a poset of given dimension from \cite{DM41}. The theme of order dimension remains to be a topic
of research interest in mathematics \cite{Wie17}.
\smallskip

A kind of opposite task to the one of generating linear extensions of a given poset $P$ on~$N$
was considered by Heath and Nema \cite{HN13} in their 2013 paper. The task considered by them
was the so-called {\em poset cover problem}, which can be re-formulated as follows.
Given a set $S$ of total orders on $N$ and a parameter $k\in {\dv N}$, decide whether
there exist posets $P_{1},\ldots, P_{k}$ on $N$ such that
$S$ is the union $\bigcup_{i} {\cal L}(P_{i})$. Of course, in case $k=1$ this reduces to the task
to decide whether there exists a poset $P$ on $N$ with $S={\cal L}(P)$.
They offered an algorithm to construct a poset~$P$
with ${\cal L}(P)\subseteq S$, where ${\cal L}(P)$ contains a prescribed total order from $S$, and
derived a result about the computational complexity of the poset cover problem.
The methodological tools they used in their paper inspired the proof of the main result in the present paper.
\medskip

The central concept in the present paper is that of the permutohedral graph over $N$,
whose nodes are enumerations of $N$. This graph is connected  and, therefore, one
can interpret it as a (finite) metric space, where the distance of two
enumerations $\permA$ and $\permB$ is the shared length of the paths between $\permA$ and $\permB$
of the least possible length, which paths are called {\em geodesics}. One can thus
consider the concept of a {\em geodetically convex set\/}, widely studied in graph-theoretical
literature~\cite{Pel13}. It is defined as such a set of nodes in the graph which, with
every two nodes, contains all the nodes on every geodesic between those two nodes.
Our main result says that a non-empty set of nodes in the permutohedral graph represents the set ${\cal L}(P)$
of linear extensions of a poset $P$ on $N$ if and only if it is a geodetically convex set.

It looks like we are not the first researchers who came with an observation of this sort. According to
Bj\"{o}rner and Wachs \cite[\S\,6]{BW91}, already in the 1970, Tits \cite[Theorem~2.19]{Tit74}
obtained a characterization of (geodetically) convex sets in the set of permutations
of an ordered set of integers $[n]:=\{1,\ldots,n\}$. To convey this result to non-specialists in Coxeter groups,
Bj\"{o}rner and Wachs came in \cite[Propositon~6.2]{BW91} with their reformulation and proof of that result by Tits.
Then, in their \cite[Theorem~6.6]{BW91}, they related these (geodetically) convex sets
to the so-called {\em labeled posets\/} on an $n$-element set, by which they meant posets on a set with a prescribed total order. Nevertheless, their presentation still remains within the framework of group theory,
because their arguments rely on some additional specialized results on the so-called weak (Bruhat) order
on the symmetric group $\sym_{[n]}$ of all permutations of the set $[n]$.

In the present paper we offer a straightforward elementary proof of the fact that linear extensions of finite posets
are characterized by geodetic convexity. Our self-contained proof does not need any specialized result on the symmetric group and is, thus, conceptually simpler. Indeed, realize that the collection of total orders
on a general finite set $N$ has no distinguished element and, therefore, it does not form
a group. Nonetheless, the mental transition to labeled posets and the relation to the results of \cite[\S\,6]{BW91} are also commented (Section~\ref{ssec.poset-encoding}).
\medskip

The permutohedral graph can be viewed as a graph with colored edges:
an equivalence on the set of its edges can be defined
by introducing the so-called {\em combinatorial labels\/} on them.
Specifically, a two-element set $\{u,v\}\subseteq N$ will be interpreted as a label
on the edge between enumerations $\permA$ and $\permB$ if these only differ
in the mutual position of $u$ and $v$. The same combinatorial label on different edges has a geometric interpretation
in terms of a certain polytope $\perN\subseteq {\dv R}^{N}$, named the {\em permutohedron}.
This is because the edges correspond to one-dimensional faces of\/ $\perN$ and
a shared label on different edges means that the edges are parallel \cite[\S\,2.3.2]{Stu24}.

The collection of all geodetically convex sets in the permutohedral graph
forms a lattice. This lattice is graded
and its height function, sometimes named a ``rank function" \cite{Sta97}, is characterized. The height of a non-empty geodetically
convex set $S$ of enumerations is shown to be the {\em number of different edge colors\/}
within the induced subgraph for $S$.
What is perhaps surprising for the reader is that, once $|N|\geq 6$, the value of this function
for $S$ \underline{does} \underline{not} \underline{coincide} with the graphical {\em diameter\/} of~$S$.
An example is given, based on a poset of dimension 3\/; in fact, the height of ${\cal L}(P)$ for a poset $P$
on $N$ coincides with the diameter of ${\cal L}(P)$ if and only if the dimension of $P$ is at most two \cite[Theorem\,4.2]{Mas09}.

It is also shown that, if $|N|\geq 2$, the above lattice
is anti-isomorphic to the {\em poset-based lattice\/} over $N$, which is the collection of all
posets on $N$ extended by the full binary relation $N\times N$. In particular, its sub-collection
of posets $P$ on $N$ forms a meet semi-lattice, which is graded and
its height function is given by the number of compared pairs (in $P$) of distinct elements of $N$.

To offer global perspective on alternative descriptions of finite posets,
two other views on them are discussed in the end of the paper. In fact, these
are special instances of alternative anti-isomorphic views on {\em preposets\/} on $N$, which are
transitive and reflexive binary relations on~$N$. The geometric approach
is to associate a poset on $N$ with a particular full-dimensional polyhedral
cone in ${\dv R}^{N}$, named the {\em braid cone\/}; this relates to \cite[\S\,3.4]{PRW08}.
The other interpretation is combinatorial, in terms of certain rings of subsets of $N$,
we name (finite) {\em topologies\/}; this relates to well-known Birkhoff's representation theorem for finite distributive lattices \cite[\S\,III.3]{Bir95}. Both interpretations are presented by means of (order-reversing) Galois connections.
\medskip

The structure of the paper is as follows.  In Section~\ref{sec.tutorial} particular
lattice-theoretical tools are recalled for the reader's convenience, while, in Section~\ref{sec.preli-other},
(most of the) specific concepts of the present paper are introduced.
In Section~\ref{sec.linear-extensions} we establish formally the correspondence between posets on $N$ and
their linear extensions. The main results are presented in Section~\ref{sec.main-results}.
Section~\ref{sec.other-views} contains remarks on two other cryptomorphic views on finite posets and conclusions.
The proofs of most results were moved to the Appendix to make the paper smoothly readable.

\section{Lattice-theoretical tutorial}\label{sec.tutorial}
As mentioned in the introduction, we are going to use some tools from lattice theory, in particular the Galois connections methodology, to present alternative views on finite posets. This means that the reader is assumed to be familiar with basic concepts from lattice theory, which can be found in relevant monographs \cite{Bir95,Sta97,Zie95}.
Nonetheless, being aware of the fact that a general reader, to which this paper is addressed, may not be acquaint with some of these tools, we have decided to recall these notions and results in this first preliminary section.
Experts in lattice theory may skip this section and consult it only if they are not sure about our terminology.

\subsection{Posets, totally ordered sets, preposets, and lattices}\label{ssec.posets-defintions}
A partially ordered set $(\latt , \preceq)$, briefly a {\em poset}, is a
non-empty set $\latt$ endowed with a partial order (= ordering), that is,
a binary relation $\preceq$ on $\latt$ which is reflexive, anti-symmetric, and transitive.
The symbol $l\prec u$ for $l,u\in \latt$ will denote $l\preceq u$ with $l\neq u$.
A {\em total order\/} is a partial order in which every two elements are
comparable: $\forall\, l,u\in \latt$ either $u\preceq l$ or $l\preceq u$.

A {\em preposet}, also {\em preorder\/} or {\em quasi-order}, is a non-empty set $\latt$ endowed with a reflexive and transitive relation $R\subseteq \latt\times\latt$.
The {\em opposite\/} of $R$ is then $R^{op}:=\{\, (u,l)\,:\, (l,u)\in R\}$, being also a preposet.
A preposet which is additionally symmetric is an {\em equivalence\/} (on $\latt$). An elementary
fact is that, given a preposet $R$, one can consider its ``inner" equivalence $E:=R\cap R^{op}$ and interpret
$R$ as a partial order $\preceq$ on the set ${\cal E}$ of equivalence classes of $E$, defined by
$$
A\preceq B\quad \mbox{for $A,B\in {\cal E}$}
\quad:=\quad
[\,\exists\, l\in A~~\exists\, u\in B \,:\ (l,u)\in R\,]\,.
$$
Indeed, then $(l,u)\in R$ if and only if $A\preceq B$ for $A,B\in {\cal E}$ determined by $l\in A$ and $u\in B$.

A partially ordered set $(\latt ,\preceq)$ is called a {\em lattice\/} \cite[\S\,I.4]{Bir95} if every two-element subset of $\latt$ has both a least upper bound, also called the {\em join\/} or the ``supremum", and a greatest lower bound, also named the {\em meet} or the ``infimum".
A finite lattice, that is, a lattice $(\latt ,\preceq)$ with $\latt$ finite, is necessarily {\em complete}, which means that the requirement above holds for any subset of\/ $\latt$, instead of a two-element one. Every complete lattice has a {\em least element},
denoted by $\mbox{\bf 0}$, which satisfies $\mbox{\bf 0}\preceq u$ for any $u\in \latt$ and a {\em greatest element}, an element $\mbox{\bf 1}\in \latt$ satisfying $l\preceq \mbox{\bf 1}$ for any $l\in \latt$.
A~{\em meet semi-lattice\/} is a poset in which a meet exists for every two-element subset.

A poset $(\latt ,\preceq)$ is (order) {\em isomorphic\/} to a poset $(\latt^{\prime},\preceq^{\prime})$
if there is a bijective mapping $\iota$ from $\latt$ onto $\latt^{\prime}$ which preserves the order,
that is, $l\preceq u$ if and only if $\iota(l)\preceq^{\prime}\iota(u)$ for $l,u\in \latt$. If this is the case
for a lattice $(\latt ,\preceq)$
then both suprema and infima are preserved by the isomorphism $\iota$ and the lattices
are considered to be identical mathematical structures.

On the other hand, $(\latt ,\preceq)$ is {\em anti-isomorphic\/} to a poset
$(\latt^{\prime},\preceq^{\prime})$ if there is a bijective mapping $\varrho$ from $\latt$ onto $\latt^{\prime}$ which reverses the order: for $l,u\in \latt$, one has $l\preceq u$ if and only if $\varrho(u)\preceq^{\prime}\varrho (l)$.
In case of a lattice, the suprema are transformed to infima and conversely.
We say then that the posets/lattices $(\latt ,\preceq)$ and $(\latt^{\prime},\preceq^{\prime})$
are (order) {\em dual\/} to each other.

\subsection{Covering relation in a poset and graded finite lattices}\label{ssec.covering-relation}
Given a poset $(\latt ,\preceq)$ and $l,u\in \latt$ such that $l\prec u$ and there is no $e\in \latt$ with
$l\prec e\prec u$, we will write $l\cover u$ and say that $l$ is {\em covered\/} by $u$ or, alternatively, that
$u$ {\em covers\/} $l$.

The covering relation is the basis for a common method to visualize finite posets, that is, posets $(\latt ,\preceq)$ with $\latt$ finite. In the so-called {\em Hasse diagram\/} of a finite poset
$(\latt ,\preceq)$, elements of $\latt$ are represented by circles, possibly with identifiers of
the elements attached. If $l\prec u$ then the circle for $u$ is placed higher than the circle for $l$.
If, moreover, $l\cover u$ then a segment is drawn between the circles.
Clearly, any finite poset is defined up to isomorphism by its Hasse diagram.

The {\em atoms\/} of a complete lattice $(\latt ,\preceq)$ are its elements that cover the least element in the lattice, that is, $a\in \latt$ satisfying $\mbox{\bf 0}\cover a$. Analogously, {\em coatoms\/} are elements of the lattice covered by the greatest element, that is, elements $c\in \latt$ satisfying $c\cover \mbox{\bf 1}$.
We will call a complete lattice $(\latt ,\preceq)$ {\em atomistic\/} if every element of
$\latt$ is the join of a set of atoms (in $\latt$). Note that some authors
use the term ``atomic" instead of ``atomistic", see \cite[\S\,VIII.9]{Bir95} or \cite[\S\,2.2]{Zie95},
but other authors use the term ``atomic" to name a weaker condition on $(\latt,\preceq)$, namely
that, for every non-least element $e\neq \mbox{\bf 0}$ of the lattice
$\latt$, an atom $a$ exists such that $a\preceq e$.
Analogously, a complete lattice $(\latt ,\preceq)$
is {\em coatomistic\/} if every element of $\latt$ is the meet of a set of coatoms.
It is evident that a lattice anti-isomorphic to an atomistic lattice is coatomistic and conversely.

A poset $(\latt ,\preceq)$ is called {\em graded\/} \cite[\S\,I.3]{Bir95}
if there exists a function $h:\latt\to {\dv Z}$, called a {\em height function},
which satisfies, for any $l,u\in \latt$,
\begin{itemize}
\item if $l\prec u$ then $h(l)<h(u)$,
\item if $l\cover u$ then $h(u)=h(l)+1$.
\end{itemize}
Note that, in the case of a finite poset, the former condition is superfluous because it follows from the latter one,
which essentially means that the circles in the Hasse diagram can be arranged into distinct layers.
A standard sufficient condition for a finite poset to be graded is in terms of {\em chains\/}, which are (non-empty)  subsets of $\latt$ that are totally ordered; the {\em length\/} of such a chain is the number
of its elements minus one.

The reader can readily see that any finite poset in which all its maximal chains (in the sense of inclusion) have the same length is graded \cite[\S\,3.1]{Sta97}.\footnote{Stanley in \cite[p.\,99]{Sta97}
defines a (finite) graded poset as a poset in which all maximal chains have the same length.
This is a strictly stronger condition than our definition taken from Birkhoff \cite[p.\,5]{Bir95}; see Example~\ref{exa.poset-diagrams}.}
It is then immediate to observe a well-known fact is that a finite lattice is graded if and only if
all maximal chains in it have the same length. Clearly, the dual lattice to a graded finite lattice is a graded lattice.

\subsection{Galois connections for representation of a complete lattice}\label{ssec.app.Galois}
There is a general method for generating examples of complete lattices
in which method a lattice together with its dual lattice is defined on the basis
of a given binary relation between two sets. The method is universal
in the sense that every complete lattice is isomorphic to a lattice
defined in this way. Moreover, every finite lattice is isomorphic to
a lattice defined on the basis of a given binary relation between elements of two finite sets.

This method uses the so-called {\em Galois connections}, also called {\em polarities} \cite[\S\,V.7]{Bir95},
and allows one to represent every finite lattice in the form of a {\em concept lattice}, whose notion has an elegant
interpretation from the point of view of formal ontology \cite{GW99}.
To present the method we first need to recall a couple of relevant notions.

Given a set $X$, its power set will be denoted by ${\cal P}(X):=\{\, S\,:\ S\subseteq X\,\}$.
A {\em Moore family\/} of subsets of $X$, also called a {\em closure system\/} (for $X$), is a collection ${\calF}\subseteq {\cal P}(X)$ of subsets of $X$ which is closed under (arbitrary) set intersection and includes the set $X$ itself\/: $\calD\subseteq {\calF} ~\Rightarrow ~ \bigcap\calD\in {\calF}$ with a convention $\bigcap\calD=X$ for $\calD=\emptyset$. An elementary observation is that any
Moore family $({\calF},\subseteq)$, equipped with the set inclusion relation as a partial order,
is a complete lattice \cite[Theorem\,2 in \S\,V.1]{Bir95}.

A {\em closure operation\/} on subsets of $X$ is
a mapping $\cl :{\cal P}(X)\to {\cal P}(X)$ which is {\em extensive}, which means $S\subseteq \cl (S)$ for $S\subseteq X$, {\em isotone}, which means
$S\subseteq T\subseteq X \Rightarrow \cl (S)\subseteq \cl (T)$, and {\em idempotent}, which means $\cl (\cl (S))=\cl (S)$ for $S\subseteq X$. A set $S\subseteq X$ is then called {\em closed with respect to $\cl$} if $S=\cl (S)$.
The basic fact is that the collection of subsets of $X$ closed
with respect to $\cl$ is a Moore family and any Moore family can be defined in this way \cite[Theorem\,1 in \S\,V.1]{Bir95}. Specifically, given a Moore family $\calF$ of subsets of $X$, the formula
$$
\cl_{{\calF}}(S) \,:=\, \bigcap\, \{\,F\subseteq X\,:\, S\subseteq F\in
{\calF}\,\} \quad \mbox{\rm for} ~ S\subseteq X\,,
$$
defines a closure operation on subsets of $X$ having ${\calF}$ as
the collection of closed sets with respect to $\cl_{{\calF}}$;
see \cite[Theorem\,1]{GW99}.

In other words, there is a one-to-one correspondence between closure operations on subsets of $X$ and Moore families of subsets of $X$, which are
examples of complete lattices relative to set inclusion relation $\subseteq$.
Thus, any closure operation $\cl$ on (subsets of) $X$ yields an example of
a complete lattice $({\calF},\subseteq)$, where ${\calF}$ is the collection
of sets closed with respect to $\cl$.
\smallskip

We now present the method of Galois connections itself. Let $X$ and $Y$ be
sets and $\mathfrak{r}$ a binary relation between elements of $X$ and $Y$: $\mathfrak{r}\subseteq X\times Y$.
We are going to write $x\,\mathfrak{r}\,y$ to denote that $x\in X$ and $y\in Y$ are in this
{\em incidence\/} relation: $(x,y)\in \mathfrak{r}$. Then every subset $S$ of $X$ can be assigned
a subset of $Y$ as follows (the {\em forward\/} direction):
$$
S\subseteq X \quad\mapsto\quad S^{\rGalo} \,:=\,  \{\, y\in Y\,:\ x\,\mathfrak{r}\,y \quad
\mbox{for every $x\in S$}\,\}~ \subseteq Y\,.
$$
Analogously, every subset $T$ of\/ $Y$ can be assigned a subset of $X$ (the {\em backward\/} direction):
$$
T\subseteq Y \quad\mapsto\quad T^{\lGalo} \,:=\,  \{\, x\in X\,:\ x\,\mathfrak{r}\,y \quad
\mbox{for every $y\in T$}\,\}~ \subseteq X\,.
$$
The mappings $S\mapsto S^{\rGalo}$ and $T\mapsto T^{\lGalo}$ are then
{\em Galois connections\/} between the lattices $({\cal P}(X),\subseteq)$ and $({\cal P}(Y),\subseteq)$.
These are {\em order-reversing\/} connections, because
$S_{1}\subseteq S_{2}\subseteq X \,\Rightarrow\, S_{1}^{\rGalo}\supseteq S_{2}^{\rGalo}$ and
$T_{1}\subseteq T_{2}\subseteq Y \,\Rightarrow\,  T_{1}^{\lGalo}\supseteq T_{2}^{\lGalo}$.
The point is that the mapping $S\mapsto
S^{\rGalo\lGalo}:=(S^{\rGalo})^{\lGalo}$ is a closure operation on subsets of $X$ and the mapping $T\mapsto
T^{\lGalo\rGalo}:=(T^{\lGalo})^{\rGalo}$ is a closure operation on subsets of $Y$.\footnote{A useful observation
is that $S^{\rGalo}=S^{\rGalo\lGalo\,\rGalo}$ for any $S\subseteq X$.}
Thus, they define a complete lattice
$({\calX}^{\mathfrak{r}},\subseteq)$ of closed subsets of $X$ and a complete lattice $({\calY}^{\mathfrak{r}},\subseteq)$ of
closed subsets of $Y$. In addition to that,
\begin{itemize}
\item one has $S\in {\calX}^{\mathfrak{r}}$ if and only if
$S=T^{\lGalo}$ for some $T\subseteq Y$, that is,\\
${\calX}^{\mathfrak{r}}$ consists of subsets of $X$ that are images of the {\em backward\/} Galois connection $\lGalo$,
\item analogously $T\in {\calY}^{\mathfrak{r}}$ if and only if $T=S^{\rGalo}$ for some $S\subseteq X$, that is,\\
${\calY}^{\mathfrak{r}}$
consists of subsets of $Y$ that are images of the {\em forward\/} Galois connection $\rGalo$, and
\item the inverse of the mapping $S\in {\calX}^{\mathfrak{r}}\mapsto S^{\rGalo}$ is the mapping
$T\in {\calY}^{\mathfrak{r}}\mapsto T^{\lGalo}$, defining an anti-isomorphism between the lattices
$({\calX}^{\mathfrak{r}},\subseteq)$ and $({\calY}^{\mathfrak{r}},\subseteq)$ \cite[Theorem\,19 in \S\,V.7]{Bir95}.\\
This anti-isomorphism is sometimes called the (respective) {\em Galois correspondence}.
\end{itemize}
In particular, the complete lattice $({\calX}^{\mathfrak{r}},\subseteq)$ and
its dual lattice $({\calY}^{\mathfrak{r}},\subseteq)$ have together been defined on the basis of a given binary relation
$\mathfrak{r}\subseteq X\times Y$. Note that in the case of a finite $X$, or of a finite $Y$, the lattice
$({\calX}^{\mathfrak{r}},\subseteq)$ is finite, too.

The arguments why this method is universal for representing complete lattices are based on the idea from \cite[\S\,V.9]{Bir95}. Specifically, given a complete lattice $(\latt ,\preceq)$, one can put $X:=\latt$,
$Y:=\latt$ and consider a binary relation $\mathfrak{r}:= \{\, (x,y)\in X\times Y\,:\  x\preceq y\}$.
Then, for any $S\subseteq X=\latt$, $S^{\rGalo}$ is the set of upper bounds for $S$ in $(\latt ,\preceq)$
while, for any $T\subseteq Y=\latt$, $T^{\lGalo}$ is the set of lower bounds for $T$.
The reader can readily see that ${\calX}^{\mathfrak{r}}$ coincides with the set of principal ideals in $(\latt ,\preceq)$, which are sets of the form $\{z\in \latt\,:\ z\preceq a\,\}$ for some $a\in \latt$.
Indeed, $S\in {\calX}^{\mathfrak{r}}$ means $S=S^{\rGalo\lGalo}$ and one can take in place of $a$
the join of $S$: then $z\preceq a$ implies $z\in S^{\rGalo\lGalo}=S$.
It is immediate that any poset $(\latt ,\preceq)$ is isomorphic with the poset of its principal ideals ordered by inclusion. Note that in case of a finite lattice $(\latt ,\preceq)$, the sets $X$ and $Y$ are both finite by the construction.

As a by-product of the above construction we observe that any complete lattice is isomorphic
to a Moore family of subsets of a non-empty set $X$ (ordered by set inclusion), and, therefore,
it can be introduced by means of a closure operation on subsets of $X$.

\section{Additional preliminaries}\label{sec.preli-other}
We also assume the reader is familiar with elementary graph-theoretical concepts and has mild
knowledge of polyhedral geometry \cite{Zie95}. Despite this, we recall our basic concepts in this section.

\subsection{Basic notational conventions}\label{ssec.notation-convention}
Throughout the (rest of the) paper, the symbol $N$ will denote a finite non-empty {\em ground set\/} of variables, while the number of its elements will be denoted by $n:=|N|\geq 1$. We intentionally regard the ground set $N$ as an \underline{unordered} labeled set, that is, a set with distinguishable elements but without any prescribed order. In other words, we interpret $N$ as a cardinal number.

Contrary to that, given a natural number $n\in {\dv N}$, the symbol $[n]\,:=\,\{1,\ldots, n\}$ will denote a totally \underline{ordered} set of integers between $1$ and $n$. So, we interpret $[n]$ as an ordinal number. The subtle reasons for these conventions, not evident at first sight, are explained in Appendix~\ref{app.remark-conventions}.
\smallskip

Recall that the power set of $N$ is denoted by $\caP := \{ A: A\subseteq N\}$.
The {\em diagonal\/} in the Cartesian product $N\times N$ 
will be denoted by $\diag:=\{\, (u,u)\,:\ u\in N\}$.
Given $u\in N$, the symbol $u$ will also be used to denote the one-element set (= singleton) $\{u\}$.
Given $A,B\subseteq N$, the symbol $AB$ will be a shortened notation for the union $A\cup B$.

The symbol $X^{Y}$, where $X,Y$ are non-empty sets, has a standard meaning in mathematics:
it denotes the set of functions from $Y$ to $X$. In particular, the
symbol ${\dv R}^{N}$ denotes the set of functions from $N$
into the set ${\dv R}$ of reals, which functions we interpret as vectors with components indexed
by the elements of $N$. Thus, given $x\in {\dv R}^{N}$ and $u\in N$, the symbol
$x_{u}$ will denote the real value assigned to $u$ by $x$, that is, the {\em component\/} of $x$ indexed by $u$.
The {\em incidence vector\/} $\chi_{A}$ for a subset $A\subseteq N$,
is a vector $x\in {\dv R}^{N}$ specified by $x_{u}:=+1$ if $u\in A$
and $x_{u}:=0$ if $u\in N\setminus A$.

Our (scarcely used) infix notation for composition of functions $\trans$ and $\permG$ uses a dot as a symbol:
$(\trans\cdot \permG)(y):= \trans(\permG (y))$ for any argument $y$ of $\permG$.

\subsection{Transitive closure and acyclic binary relations}\label{ssec.trans-closure}
Because the intersection of transitive relations on $N$ is transitive, for every
binary relation $R\subseteq N\times N$, there exists the least transitive relation $T\subseteq N\times N$ with $R\subseteq T$, which is called the {\em transitive closure of $R$} and denoted by $\tr (R)$. It is easy to see that
$(u,v)\in \tr(R)$ if and only if there exists a sequence $u=u_{1},\ldots, u_{k}=v$, $k\geq 2$, in $N$ with
$(u_{i},u_{i+1})\in R$ for $i=1,\ldots, k-1$; additionally, without loss of generality
one can assume that both the sequence $u_{1},\dots,u_{k-1}$ and the sequence $u_{2},\dots,u_{k}$ consists of pairwise distinct
elements. Thus, in case $u\neq v$, this implies that
$u_{1},\dots,u_{k}$ could be chosen to have pairwise distinct elements.

We say that a binary relation $R\subseteq N\times N$ is {\em acyclic\/} if there is no sequence $u_{1},u_{2},\ldots, u_{k}$, $k\geq 2$, of elements of $N$ such that $u_{k}=u_{1}$ and $(u_{i},u_{i+1})\in R$ for $i=1,\ldots k-1$. This basically means that $\diag\cap\tr(R)=\emptyset$.
A well-known equivalent definition of acyclicity is that all elements of $N$ can be ordered into a sequence $u_{1},\ldots,u_{n}$, $n=|N|$, such that $(u_{i},u_{j})\in R$ implies $i<j$. Indeed,
the necessity can be verified by induction on $n$, because every acyclic relation $R$ has a terminal variable,
which is $t\in N$ such that there is no $w\in N$ with $(t,w)\in R$.

For a transitive relation $T\subseteq N\times N$, the next two conditions are equivalent to acyclicity:
\begin{itemize}
\item $T$ is {\em irreflexive}, which means that there is no $u\in N$ with $(u,u)\in T$,
\item $T$ is {\em asymmetric}, which means that $(u,v)\in T$ implies $(v,u)\not\in T$.
\end{itemize}
Indeed, realize that, for a transitive relation $T$, negations of those conditions are equivalent.\\[0.3ex]
These observations allow one to verify easily the following characterization of posets on $N$.

\begin{fact}\label{fct.posets}\rm
A relation $P\subseteq N\times N$ is a poset if and only if $\diag\subseteq P$ and
the strict version $P\setminus\diag$ of this relation is both transitive and acyclic.
\end{fact}
\smallskip

Given a poset $P$ on $N$ and $(u,v)\in N\times N$, we will sometimes apply the next symbols:
\begin{itemize}
\item $u\preceq v\,\,[P]$ for {\em non-strict comparison}, that is, $(u,v)\in P$,
\item $u\prec v\,\,[P]$ for {\em strict comparison}, that is, $(u,v)\in P\setminus\diag$,
\item $u\cover v\,\,[P]$ for {\em covering}, that is, $(u,v)\in P\setminus\diag$, $\neg[\,\exists\, w\in N\setminus \{u,v\} ~:~ (u,w),(w,v)\in P\,]$,
\item $u\|v\,[P]$ for {\em incomparability\/} of $u$ and $v$ in $P$, that is,
$(u,v)\not\in P$ and $(v,u)\not\in P$.
\end{itemize}
When the poset $P$ is clear from the context, the symbol $[P]$ might be omitted.

\subsection{Finite posets depicted as directed acyclic graphs}\label{ssec.poset-as-DAG}
There are two (different) standard  ways to represent partial orders on a non-empty finite ground
set $N$ by means of directed acyclic graphs having $N$ as the set of nodes.
The point is that every acyclic binary relation $R$ on $N$ can be depicted in the form of a {\em directed acyclic graph\/} having $N$ as the set of nodes, where an {\em arrow\/} (= directed edge) $u\to v$ is drawn if and only if $(u,v)\in R$.

The first way to represent a poset $P$ on $N$ is based on its {\em strict version}, that is, $R:=P\setminus\diag$.
Owing to the characterization of posets on $N$ in Fact~\ref{fct.posets}, it leads
to a {\em transitive directed acyclic\/} graph with node set $N$. Formulated in other words,
transitive directed acyclic graphs over $N$ are in obvious one-to-one correspondence with posets on $N$.

The second way comes from the {\em covering relation\/} $C:= \{\, (u,v)\in N\times N\,:\ u\cover v\,\,[P]\,\}$. Since $C\subseteq P\setminus\diag$, it is also acyclic.
The respective directed graph is (strongly) {\em anti-transitive\/} in the sense that $(u_{i},u_{i+1})\in C$
for $i=1,\ldots k-1$, $k\geq 3$, implies $(u_{1},u_{k})\not\in C$. This yields the so-called
{\em directed Hasse diagram\/} of the poset $P$, whose underlying undirected graph is called the {\em cover graph\/} of $P$ in \cite{Wie17}.
The following elementary observation is left to the reader.

\begin{fact}\label{fct.covering}\rm
The covering relation $C$ of a finite poset $P$ is the least subset of $P\setminus\diag$ with $\tr(C)=P\setminus\diag$.
\end{fact}

This is to illustrate the pictorial ways of representing finite posets.

\begin{example}\label{exa.poset-diagrams}\rm
Let us put $N:=\{a,b,c,d\}$ and
$P:= \diag\cup \{\, (a,b), (a,d), (b,d),(c,d)\,\}$. Then $P$ is clearly a poset on $N$.
In Figure~\ref{fig.poset-diagrams}(a) the respective transitive directed acyclic graph is shown.
For reader's convenience, the covering pairs are shown in bold there. Figure~\ref{fig.poset-diagrams}(b) then contains the respective directed Hasse diagram, obtained by the removal of the only non-covering pair $(a,d)$ in Figure~\ref{fig.poset-diagrams}(a).
To relate this to the customary conventional depicting of Hasse diagrams, as presented in \cite[p.\ 98]{Sta97},
we show in Figure~\ref{fig.poset-diagrams}(c) the corresponding classic Hasse diagram, which is undirected.
By directing the graph in Figure~\ref{fig.poset-diagrams}(c) from bottom to top one gets
what is in Figure~\ref{fig.poset-diagrams}(b), except for the placement of the nodes.
The technical difference is perhaps that rotating a directed Hasse diagram does not result in a loss of information.
Note that the poset $P$ discussed in this example is graded according to our definition from Section~\ref{ssec.covering-relation} despite it contains two maximal chains of different length, namely $c\cover d$ and $a\cover b\cover d$.
\end{example}

\begin{figure}[t]
\setlength{\unitlength}{0.8mm}
\begin{center}
\begin{picture}(40,40)
\put(7,20){\circle{3}}
\put(20,7){\circle{3}}
\put(20,33){\circle{3}}
\put(33,20){\circle{3}}
%
\put(18.5,31.5){\vector(-1,-1){10}}
%
\put(20,37){\makebox(0,0){$a$}}
\put(2,20){\makebox(0,0){$d$}}
\put(37,20){\makebox(0,0){$b$}}
\put(20,3){\makebox(0,0){$c$}}
\put(4,36){\makebox(0,0){\bf (a)}}
\thicklines
\put(18.5,8.5){\vector(-1,1){10}}
\put(21.5,31.5){\vector(1,-1){10}}
\put(31,20){\vector(-1,0){22}}

\end{picture}
\hspace*{5mm}
\begin{picture}(40,40)
\put(7,20){\circle{3}}
\put(20,7){\circle{3}}
\put(20,33){\circle{3}}
\put(33,20){\circle{3}}
\put(18.5,8.5){\vector(-1,1){10}}
\put(21.5,31.5){\vector(1,-1){10}}
\put(31,20){\vector(-1,0){22}}
\put(20,37){\makebox(0,0){$a$}}
\put(2,20){\makebox(0,0){$d$}}
\put(37,20){\makebox(0,0){$b$}}
\put(20,3){\makebox(0,0){$c$}}
\put(4,36){\makebox(0,0){\bf (b)}}
\end{picture}
\hspace*{5mm}
\begin{picture}(53,40)
\put(7,20){\circle{3}}
\put(20,33){\circle{3}}
\put(33,20){\circle{3}}
\put(46,7){\circle{3}}
\put(18.5,31.5){\line(-1,-1){10}}
\put(21.5,31.5){\line(1,-1){10}}
\put(34.5,18.5){\line(1,-1){10}}
\put(20,37){\makebox(0,0){$d$}}
\put(3,20){\makebox(0,0){$c$}}
\put(37,21){\makebox(0,0){$b$}}
\put(50,7){\makebox(0,0){$a$}}
\put(4,36){\makebox(0,0){\bf (c)}}
\end{picture}
\caption{Three ways to graphically represent the same poset $P$ over $N:=\{a,b,c,d\}$.\label{fig.poset-diagrams}}
\end{center}
\end{figure}

\subsection{Enumerations and permutohedron}\label{ssec.enum-permutohedron}
Ordered lists of (all) elements of $N$ (without repetition) will be called {\em enumerations\/} (of $N$).
Formally, an enumeration of $N$ is a one-to-one (= bijective) mapping $\permA : [n]\to N$ from $[n]$ onto~$N$.
We will use the record $|\,\permA(1)|\ldots|\permA(n)\,|$ to specify such an enumeration.
Every enumeration of $N$ can be interpreted as a total order on $N$. In fact, there are two standard ways to do that:
one can either interpret $\permA$ in the ascending way as $\permA(1)\prec \ldots\prec \permA(n)$
or one may interpret it in the descending way as $\permA(1)\succ \ldots\succ \permA(n)$.

The set of (all) enumerations of $N$ will be denoted by $\enuN$.
Every enumeration corresponds uniquely to the respective {\em rank vector\/}, which is
a vector in $x\in {\dv R}^{N}$ defined by $x_{u}:= \permA_{-1}(u)$ for $u\in N$.
In other words, the rank vector (for an enumeration $\permA$)
is the inverse $\permA_{-1}: N\to [n]$ of the mapping $\permA : [n]\to N$. We are going
to denote it by $\rankv_{\permA}$ in a geometric context (of ${\dv R}^{N}$).

The {\em permutohedron\/} in ${\dv R}^{N}$, denoted by $\perN$, is defined as the convex hull
of the rank vectors for (all) enumerations of $N$. Formally,
$\perN :=\conv (\,\{\rankv_{\permA}\,:\ \permA\in\enuN\}\,)$. It is a subset of
the affine space $\{\, y\in {\dv R}^{N}\,:\ \sum_{u\in N} y_{u}={n+1\choose 2}\,\}$; therefore,
its (Euclidean) dimension is $n-1$.

\subsection{Face lattice of the permutohedron}\label{ssec.face-lattice}
Recall that the set of (all) faces of any polytope, ordered by inclusion, is a finite lattice, known to be both atomistic and coatomistic. In addition to that, this {\em face lattice\/} is a graded lattice, where the dimension (of faces) plays the role of a height function \cite[Theorem\,2.7]{Zie95}.

In this subsection we present an elegant {\em combinatorial\/} description of (non-empty) faces of
the permutohedron $\perN$, mentioned in \cite[p.\,18]{Zie95}, whose proof can be found in \cite[\S\,1]{BS96}.
Every non-empty face of the permutohedron $\perN$ corresponds to an ordered partition
of the ground set~$N$ into non-empty blocks. We will use the symbol $|\,A_{1}|\ldots|A_{m}\,|$ to
denote such a partition into $m$ blocks, $1\leq m\leq n$. Note that this notation is consistent
with our notation for enumerations. In fact, the ordered partition serves as a compact
description of the set of those enumerations of $N$ which encode the {\em vertices\/} of the
corresponding face of $\perN$: these are the respective rank vectors.
Specifically, it is the set of those enumerations $\permA$ of $N$ which are consistent with $|\,A_{1}|\ldots|A_{m}\,|$: these are the ordered lists $\permA$ of elements of $N$ where first (all)
the elements of $A_{1}$ are listed, then (all) the elements of $A_{2}$ follow and so on.
Thus, any face of\/ $\perN$ also corresponds to a special subset of\/ $\enuN$.
The (Euclidean) dimension of the face described by $|\,A_{1}|\ldots|A_{m}\,|$ is $n-m$ then.
Later, in Section~\ref{sec.linear-extensions}, we interpret such (non-empty) face-associated subsets of $\enuN$
as a special case of subsets of\/ $\enuN$ assigned to partial orders on $N$.

\subsection{Some graph-theoretical concepts}
We assume that the reader is familiar with elementary graph-theoretical notions and
understands the next definitions, which apply to any {\em connected undirected graph\/} with
a node set~$U$.

A {\em geodesic\/} between nodes $\permA,\permB\in U$ is a walk between $\permA$ and $\permB$ (in the graph)
which has the shortest possible length among such walks. Note that it is necessarily a path
(= nodes are not repeated in it) and that several geodesics may exist between two nodes. The (graphical)
{\em distance\/} $\dist(\permA,\permB)$ between two nodes $\permA,\permB\in U$ is the length of a geodesic between them. It makes no problem to observe that it satisfies the axioms of a metric and defines a finite metric space.

We say that a node $\permC\in U$ {\em is between nodes\/} $\permA\in U$ and $\permB\in U$ if $\permC$ belongs to some geodesic between $\permA$ and $\permB$; an equivalent condition is that $\dist(\permA,\permC)+\dist(\permC,\permB)= \dist(\permA,\permB)$.
A set $S\subseteq U$ is {\em geodetically convex\/} if, for any $\permA,\permB\in S$, all nodes between them belong to $S$. It is evident that the intersection of geodetically convex sets is a geodetically convex set.

\subsection{Permutohedral graph}\label{ssec.permut-graph}
We say that enumerations $\permA,\permB\in\enuN$ differ by an {\em adjacent transposition\/} if
\begin{eqnarray}
\lefteqn{\hspace*{-2cm}\exists\, 1\leq i<n ~~~\mbox{such that}~~ \permA(i)=\permB(i+1),~ \permA(i+1)=\permB(i),} \label{eq.adjacent-transpos}\\
&\mbox{and}&
\permA(k)=\permB(k)\qquad \mbox{for remaining $k\in [n]\setminus\{i,i+1\}$}\,. \nonumber
\end{eqnarray}
We are going to interpret the set of enumerations of $N$ as an undirected {\em permutohedral graph\/}
in which $\enuN$ is the set of nodes and edges are determined by adjacent transpositions.
The reason for this terminology is as follows. Despite the given definition of adjacency
for enumerations being fully combinatorial, it has a natural geometric interpretation.
Enumerations $\permA,\permB\in\enuN$ are adjacent in this graph if and only if the segment
$[\rankv_{\permA},\rankv_{\permB}]\subseteq {\dv R}^{N}$ between the respective rank vectors is a face
(= a geometric edge) of the permutohedron $\perN$. This follows from the
results presented in Section~\ref{ssec.face-lattice}: the geometric edges of $\perN$ are described
by ordered partitions of the form $|\,A_{1}|\ldots|A_{n-1}\,|$; then $|A_{i}|=2$ for unique $i\in[n-1]$ and
$|A_{j}|=1$ for remaining $j\neq i$.

To prove our main result (Theorem~\ref{thm.poset-characterization})
it is instrumental to introduce a particular {\em edge-labeling\/}  (= edge coloring) for this graph.
Specifically, we are going to label its edges by two-element subsets $\{u,v\}$ of
our ground set $N$: if $\permA,\permB\in\enuN$ differ by an adjacent transposition
\eqref{eq.adjacent-transpos} then we label the edge between $\permA$ and $\permB$ by the set
$\{u,v\}:=\{\permA(i),\permA(i+1)\}=\{\permB(i),\permB(i+1)\}$.
We will call the set $\{u,v\}$ a {\em combinatorial label\/}\label{df.comb-label} on the edge and denote this by
$\permA \stackrel{uv}{\longleftrightarrow} \permB$.

If specifically $u=\permA(i)$ and $v=\permA(i+1)$ then \eqref{eq.adjacent-transpos} yields
$\rankv_{\permA}-\rankv_{\permB}=\chi_{v}-\chi_{u}$, that is, the rank vectors only differ in the components
$u$ and $v$ (by one). In particular, if $\permA \stackrel{uv}{\longleftrightarrow} \permB$  and $\permC \stackrel{uv}{\longleftrightarrow} \permD$
then $\rankv_{\permA}-\rankv_{\permB}=\pm (\rankv_{\permC}-\rankv_{\permD})$ meaning that the segments
$[\rankv_{\permA},\rankv_{\permB}]$ and $[\rankv_{\permC},\rankv_{\permD}]$ are {\em parallel\/},  or worded
in a different way, they are perpendicular to the same hyperplane ${\mathsf H}_{uv}:=\{\,z\in {\dv R}^{N}\, :\ z_{u}=z_{v}\,\}$.
Thus, the respective {\em combinatorial equivalence\/} of edges in the permutohedral graph, defined
as having the same combinatorial label, has a natural geometric interpretation.

\begin{figure}[t]
\setlength{\unitlength}{0.8mm}
\begin{center}
\begin{picture}(85,40)
\put(20,20){\circle{3}}
\put(28,7){\circle{3}}
\put(28,33){\circle{3}}
\put(43,7){\circle{3}}
\put(43,33){\circle{3}}
\put(51,20){\circle{3}}
\put(11,20){\makebox(0,0){$|b|a|c|$}}
\put(25,1){\makebox(0,0){$|a|b|c|$}}
\put(25,39){\makebox(0,0){$|b|c|a|$}}
\put(46,1){\makebox(0,0){$|a|c|b|$}}
\put(46,39){\makebox(0,0){$|c|b|a|$}}
\put(60,20){\makebox(0,0){$|c|a|b|$}}
\put(29.9,7){\line(1,0){11.2}}
\put(29.9,33){\line(1,0){11.2}}
\put(20.95,21.5){\line(3,5){6}}
\put(20.95,18.5){\line(3,-5){6}}
\put(50.05,21.5){\line(-3,5){6}}
\put(50.05,18.5){\line(-3,-5){6}}
\put(35.5,4){\makebox(0,0){\scriptsize $bc$}}
\put(35.5,36){\makebox(0,0){\scriptsize $bc$}}
\put(21,12.5){\makebox(0,0){\scriptsize $ab$}}
\put(50,12.5){\makebox(0,0){\scriptsize $ac$}}
\put(21,27.5){\makebox(0,0){\scriptsize $ac$}}
\put(50,27.5){\makebox(0,0){\scriptsize $ab$}}
\end{picture}
\caption{The permutohedral graph over $N:=\{a,b,c\}$ with combinatorial labels.\label{fig.hexagon}}
\end{center}
\end{figure}

\begin{example}\label{exa.permut-graph-1}\rm
To illustrate these notions put  $N:=\{a,b,c\}$. The corresponding permutohedral graph is shown
in Figure~\ref{fig.hexagon}; in fact, it is a hexagon. The nodes have the respective enumerations attached and the edges have
the combinatorial labels. We also observe that the opposite edges in the hexagon have the same label and are, thus, combinatorially equivalent.
\end{example}

Every two-element subset $\{u,v\}$ of $N$ also defines an {\em automorphism\/} of the
permutohedral graph. Specifically, the transposition $\trans_{uv}:N\to N$ of $u$ and $v$,
defined by $\trans_{uv}(u)=v$, $\trans_{uv}(v)=u$, and $\trans_{uv}(t)=t$ for $t\in N\setminus\{u,v\}$,
is a (self-inverse) permutation of $N$ and every enumeration $\permG :[n]\to N$ can be assigned its composition
$\trans_{uv}\cdot\permG$ with $\trans_{uv}:N\to N$. Note that, thanks to our distinguishing (unordered) $N$ from
(ordered) $[n]$, there is a unique way to compose these two bijections.
The mapping $\permG\mapsto \trans_{uv}\cdot\permG$
is a self-inverse transformation of\/ $\enuN$ preserving the edges in the permutohedral graph. In fact, it maps
bijectively the set
$$
S_{u\prec v}\,:=\,\{\,\permG\in\enuN \,:\ \permG_{-1}(u)<\permG_{-1}(v)\,\}
\quad\mbox{of enumerations in which $u$ precedes $v$}
$$
to its complement
$S_{v\prec u}\,=\,\{\,\permF\in\enuN \,:\ \permF_{-1}(v)<\permF_{-1}(u)\,\}$.

\begin{example}\label{exa.permut-graph-2}\rm
To illustrate these automorphisms consider again the case $N:=\{a,b,c\}$ and the permutohedral graph
in Figure~\ref{fig.hexagon}. Then the subset $\{b,c\}$ defines a graphical automorphism of the hexagon:
$|a|b|c|\leftrightarrow |a|c|b|$, $|b|a|c|\leftrightarrow |c|a|b|$, and $|b|c|a|\leftrightarrow |c|b|a|$.
This can be viewed as a horizontal swapping of the nodes in Figure \ref{fig.hexagon}.
Note in this context that one can distinguish 12 automorphisms of this graph, meaning that its symmetry group has 12 elements.
\end{example}

\section{Linear extensions}\label{sec.linear-extensions}
In this section we define a particular {\em poset-based lattice}, extending the collection
of posets on~$N$. To explain what we mean by an extension note that, in the case $n=|N|\geq 2$, the collection of posets on $N$ is not a lattice, just a meet semi-lattice. It can, however, be turned into a lattice by adding a greatest element, namely the (cyclic) relation $N\times N$.
To demonstrate readily that this extended structure indeed forms a lattice we employ the method of Galois connections, as described in Section~\ref{ssec.app.Galois}.
\smallskip

To this end we introduce a binary relation $\preced$ between elements of the set $X:=\enuN$ of enumerations
of $N$ and the elements of the Cartesian product $Y:=N\times N$. More specifically, if $\permA :[n]\to N$ is an enumeration of $N$ and $(u,v)\in N\times N$ then we consider $\permA$ to be in the incidence relation $\preced$ with the pair $(u,v)$ if and only if {\em $u$ precedes $v$\/} in $\permA$\,:
$$
\permA\,\preced\, (u,v) \quad :=\quad \permA_{-1}(u)\leq \permA_{-1}(v)\,,
\quad \mbox{written also in the form $u\preceq_{\permA} v$.}
$$
We will denote Galois connections based on this precedence relation $\preced$
using larger triangles $\rGprec$ and $\lGprec$, in order to distinguish them
from later Galois connections (in Section~\ref{sec.other-views}).
Thus, any enumeration $\permA\in\enuN$ of $N$ is assigned (by forward Galois connection) a total order on $N$:
$$
T_{\permA} := (\{\permA\})^{\rGprec} =
\{\, (u,v)\in N\times N\, :\ u\preceq_{\permA} v\,\}\,.
$$
We also introduce special notation for the case that {\em $u$ strictly precedes $v$\/} in $\permA$\,:
$$
u\prec_{\permA} v \qquad \mbox{will denote~ $\permA_{-1}(u) < \permA_{-1}(v)$,\quad that is,~ $u\prec v\,\,[T_{\permA}]$.}
$$
Analogously, in the case that {\em $u$ immediately precedes $v$\/} in $\permA$ we will write
$$
u\covA v \qquad \mbox{to denote~ $\permA_{-1}(v)-\permA_{-1}(u)=1$,\quad that is,~ $u\cover v\,\,[T_{\permA}]$.}
$$
Any subset $T$ of\/ $Y=N\times N$ is assigned (by backward Galois connection) an enumeration set
$$
{\cal L}(T):= T^{\lGprec}= \{\, \permA\in\enuN\,:\ \forall\,(u,v)\in T~~ u\preceq_{\permA} v\,\}\,,
$$
which is interpreted as the set of {\em linear extensions\/} of $T$ (provided $T\setminus\diag$ is an acyclic
relation on~$N$, as otherwise ${\cal L}(T)$ turns out to be empty).

\begin{figure}[t]
\setlength{\unitlength}{0.8mm}
\begin{center}
\begin{picture}(40,40)
\put(7,20){\circle{3}}
\put(20,7){\circle{3}}
\put(20,33){\circle{3}}
\put(33,20){\circle{3}}
%
\put(18.5,31.5){\vector(-1,-1){10}}
%
\put(20,37){\makebox(0,0){$a$}}
\put(2,20){\makebox(0,0){$d$}}
\put(37,20){\makebox(0,0){$b$}}
\put(20,3){\makebox(0,0){$c$}}
\put(4,36){\makebox(0,0){\bf (a)}}
\thicklines
\put(18.5,8.5){\vector(-1,1){10}}
\put(21.5,31.5){\vector(1,-1){10}}
\put(31,20){\vector(-1,0){22}}
\end{picture}
\hspace*{22mm}
\begin{picture}(40,60)
\put(8,36){\makebox(0,0){\bf (b)}}
\put(10,20){\circle{3}}
\put(30,20){\circle{3}}
\put(50,20){\circle{3}}
\put(12,20){\line(1,0){16}}
\put(32,20){\line(1,0){16}}
\put(7,14){\makebox(0,0){$|a|b|c|d|$}}
\put(30,14){\makebox(0,0){$|a|c|b|d|$}}
\put(53,14){\makebox(0,0){$|c|a|b|d|$}}
\put(20,23){\makebox(0,0){\scriptsize $bc$}}
\put(40,23){\makebox(0,0){\scriptsize $ac$}}
\end{picture}
\caption{The poset from Example~\ref{exa.poset-diagrams} and (the graph of) its linear extensions.\label{fig.poset-extensions}}
\end{center}
\end{figure}

\begin{example}\label{exa.poset-extensions}\rm
To illustrate the above concepts consider again the poset $P$ from Example~\ref{exa.poset-diagrams},
whose transitive directed acyclic graph is shown in Figure~\ref{fig.poset-extensions}(a). The
set of linear extensions ${\cal L}(P)$ contains three enumerations. These are depicted in Figure~\ref{fig.poset-extensions}(b), where
additionally the edges of the permutohedral graph between them are shown, together with the respective combinatorial
labels. This is what is named the {\em linear extension graph\/} (of $P$) in \cite{Mas09}.
\end{example}

\begin{lemma}\label{lem.poset-lattice}\rm
Let $({\calX}^{\preced},\subseteq )$ and $({\calY}^{\preced},\subseteq )$ denote the finite lattices defined using the Galois connections $\lGprec$ and $\rGprec$ based on the incidence relation $\preced$ between $X$ and $Y$ introduced above.
\begin{itemize}
\item[(i)] If $n=|N|\geq 2$ then\, [\,$T\in {\calY}^{\preced}$ and\, $T\subset N\times N$\,]
if and only if\/ $T$ is a poset on $N$.
\item[(ii)] One has [\,$S\in {\calX}^{\preced}$ and $S\neq\emptyset$\,] if and only if
a (uniquely determined) poset $(N,\preceq)$ exists such that $S=\{\,\permA\in\enuN\,:\
\permA_{-1}(u)<\permA_{-1}(v)~~ \mbox{whenever $u\prec v$}\,\}$, that is,\\
$S={\cal L}(P)$ for a poset~$P$ on $N$.
\item[(iii)] Any $S\subseteq\enuN$ representing a non-empty face of\/ $\Pi (N)$
(as in Section\,\ref{ssec.face-lattice})  belongs to ${\calX}^{\preced}$.
\item[(iv)] The lattice $({\calX}^{\preced},\subseteq )$ is
atomistic, coatomistic, and graded.
\end{itemize}
Coatoms in ${\calX}^{\preced}$ are precisely the sets $S_{u\prec v}\,:=\,\{\,\permA\in\enuN \,:\ u\prec_{\permA}v\,\}$,
where $u,v\in N$, $u\neq v$.
\end{lemma}
\medskip

To make the present paper smoothly readable, our elementary proof, based solely on the observations from
Section~\ref{ssec.app.Galois}, was moved to Appendix~\ref{app.lem.poset-lattice}.

Note in this context that the class ${\calX}^{\preced}\setminus\{\emptyset\}$ of poset-interpretable sets
of enumerations  can alternatively be introduced by means of a ``strict" incidence relation $\permA\odot(u,v)$
defined by $\permA_{-1}(u)<\permA_{-1}(v)$.
Except for the degenerate case $|N|=1$, the resulting lattice $({\calX}^{\odot},\subseteq)$
coincides with $({\calX}^{\preced},\subseteq)$ while the elements of\/ ${\calY}^{\odot}$
are strict versions of posets on $N$ instead.\smallskip

The reader may ask what is the closure operation corresponding to the {\em poset-based lattice\/} $({\calY}^{\preced},\subseteq)$.
Given $T\subseteq Y=N\times N$, one can show that
$$
T^{\lGprec\rGprec} ~=~
\left\{
\begin{array}{cl}
\tr(T\cup\diag) & \mbox{if $T\setminus\diag$ is acyclic},\\
N\times N & \mbox{otherwise}.
\end{array}
\right.
$$

\section{Main results}\label{sec.main-results}
In this section we present our results on linear extensions.
Section~\ref{ssec.geodesics} provides basic observations on geodesics in the permutohedral graph.
The main theorem, being a graphical characterization of poset-interpretable enumeration sets, is presented in Section~\ref{ssec.characterization}.
Section~\ref{ssec.graphical-comparison} then claims that basic order-theoretical relations for pairs of variables in a poset $P$
can be graphically recognized based on the set ${\cal L} (P) $ of its linear extensions.
Section~\ref{ssec.poset-dimension} provides an example that the height function of the lattice of geodetically convex sets in the permutohedral graph is not given by the graphical diameter for $n\geq 6$; a relevant concept of the  {\em dimension of a poset} is discussed there, too. Section~\ref{sec.main-results} is concluded by remarks, namely on an ``algorithmic" characterization of geodetically convex sets presented in \cite{HN13} (Section~\ref{ssec.remark-alter})
and on encoding posets relative to a fixed reference total order of $N$ applied in \cite{BW91}
(Section~\ref{ssec.poset-encoding}).

\subsection{Geodesics in the permutohedral graph}\label{ssec.geodesics}
Given $\permA,\permB\in\enuN$, a two-element subset $\{u,v\}$ of $N$ will be called
an {\em inversion between\/} $\permA\in\enuN$ and $\permB\in\enuN$ if the mutual orders of $u$ and $v$
in $\permA$ and in $\permB$ differ. We will denote the set of all inversions between $\permA$ and $\permB$ as follows:
\begin{eqnarray*}
\Inv[\permA,\permB] ~:=~
\left\{\,\, \{u,v\}\subseteq N\, :~
\mbox{either $[\, u\prec_{\permA}v ~\&~ v\prec_{\permB}u\,]$ ~or~
$[\, v\prec_{\permA}u ~\&~ u\prec_{\permB}v\,]$}\,\,\right\}\,.
\end{eqnarray*}
The next lemma brings some basic auxiliary observations on the permutohedral graph,
which will serve as elementary steps in the proofs of more advanced findings.

\begin{lemma}\label{lem.permut-graph}\rm
Consider $\permA,\permB\in\enuN$. Then
\begin{itemize}
\item[(i)] $\dist(\permA,\permB)=|\,\Inv[\permA,\permB]\,|$ is the number of inversions
between $\permA$ and $\permB$. In particular, $\permA$ and~$\permB$ are adjacent
in the permutohedral graph if and only if they differ by one inversion only. The label on the edge
between $\permA$ and $\permB$ is then the only inversion between them.
\item[(ii)]
A walk in the permutohedral graph is a geodesic (between $\permA$ and $\permB$) if and only if no label
on its edges is repeated. If this is the case then the set of labels on the geodesic coincides
with the set $\Inv[\permA,\permB]$ of inversions between $\permA$ and $\permB$.
\item[(iii)]
An enumeration $\permC\in\enuN$ is between $\permA$ and $\permB$ if and only if $\Inv[\permA,\permC]\cap\Inv[\permC,\permB]=\emptyset$.
Another equivalent condition is $T_{\permA}\cap T_{\permB}\subseteq T_{\permC}$, that is,
$[\, u\prec_{\permA} v ~\&~  u\prec_{\permB} v\,] ~\Rightarrow ~ u\prec_{\permC} v$.
\end{itemize}
\end{lemma}

Note that observations of this sort happen to occur independently in the literature in various contexts.
Similar or analogous results can be found in the publications on posets and on Coxeter groups, see \cite[Theorem\,1(i)]{PRRZ95}, \cite[Lemma\,2.2]{Mas09}, and \cite[Theorem\,5.1.6]{BGW03}.\footnote{By the way,
an observant reader might notice that the term ``generalized permutohedron" in poset theory~\cite{PRRZ95} refers to something else than the same term in polyhedral geometry \cite{PRW08,MUWY18} or in Coxeter groups \cite[\S\,5.16]{BGW03}.}
One of the reviewers pointed out that the claims in Lemma \ref{lem.permut-graph} follow readily from
results well known to experts in Coxeter groups. Thus, a reader familiar with those fields can consult the specialized literature. Nevertheless, in the spirit of Occam's razor, we decided to offer the reader a conceptually simpler, elementary proof. To ensure the flow of the text, we have moved this proof to Appendix \ref{app.lem.permut-graph}.

\subsection{Characterization of poset-interpretable sets of enumerations}\label{ssec.characterization}
Two special concepts for a non-empty set\/ $\emptyset\neq S\subseteq\enuN$ of nodes in the permutohedral graph are needed in the sequel.
The set of {\em inversions within $S$\/} is as follows:
$$
\Inv (S) ~:=~ \left\{\,\, \{u,v\}\subseteq N\, :\
\mbox{there exist $\permA,\permB\in S$ such that $\permA \stackrel{uv}{\longleftrightarrow} \permB$}\,\, \right\}.
$$
Note that $\Inv (S)$ can be viewed as a kind of generalization of the formerly defined set $\Inv[\permA,\permB]$ of inversions between two enumerations $\permA$ and $\permB$: it follows from Lemma~\ref{lem.permut-graph}(ii) that if one takes the set of enumerations between $\varepsilon$ and $\eta$ in place of $S$ then $\Inv(S)=\Inv[\permA,\permB]$.

The {\em covering relation for\/} $S$ is the following binary relation on $N$:
$$
\Cov (S) \,:=\, \left\{\,\, (u,v)\in N\times N\, :\
\mbox{there are $\permA\in S$, $\permB\in \enuN\setminus S$ with $\permA \stackrel{uv}{\longleftrightarrow} \permB$ and $u\prec_{\permA}v$}\,\, \right\}.
$$
Note that the latter terminology is motivated by the fact that $\Cov (S)$ appears to coincide with the covering relation $\cover$ for the poset on $N$ assigned to $S$ (see Section~\ref{ssec.graphical-comparison}).
These two notions are substantial in the proof of the next result, which can be found in Appendix~\ref{app.thm.poset-characterization}.

\begin{thm}\label{thm.poset-characterization}\em
If $n=|N|\geq 2$ then, given $S\subseteq\enuN$, $S\in {\calX}^{\preced}$ if and only if $S$ is geodetically convex.
In particular, $S={\cal L}(P)$ for a poset $P$ on $N$ if and only if $S$ is a non-empty geodetically convex set
of enumerations of $N$.
\end{thm}

\subsection{Graphical recognition of incomparable and comparable pairs}\label{ssec.graphical-comparison}
These are the observations that were implicitly behind \cite[\S\,4]{HN13}.

\begin{thm}\label{thm.poset-to-extensions}\em
Let $P$ be a poset over $N$ and $S={\cal L}(P)$. Then, for any pair $u,v\in N$, one has
\begin{itemize}
\item[(i)] $u\|v\,[P]$ if and only if $\{u,v\}\in \Inv (S)$,
\item[(ii)] $u\cover v\,\,[P]$ if and only if  $(u,v)\in\Cov (S)$,
\item[(iii)] $(u,v)\in P\setminus\diag$ if and only if  $(u,v)\in\tr(\Cov (S))$,
\item[(iv)]  $(u,v)\in P$ if and only if  $(u,v)\in\tr(\diag\cup\Cov (S))$.
\end{itemize}
\end{thm}

Our technical proof is in Appendix~\ref{app.thm.poset-to-extensions}.
It has the following easy consequence.

\begin{corol}\label{cor.exclusive}\em
Given a poset $P$ on $N$ and  $S={\cal L}(P)$, any $(u,v)\in N\times N\setminus\diag$
satisfies {\em exclusively} one of the following three conditions:
\begin{itemize}
\item[(a)] $\{u,v\}\in \Inv (S)$,
\item[(b)] $(u,v)\in\tr(\Cov (S))$,
\item[(c)] $(v,u)\in\tr(\Cov (S))$.
\end{itemize}
\end{corol}

\begin{proof}
By claims (i) and (iii) in Theorem~\ref{thm.poset-to-extensions}, the conditions are
equivalent to $u\|v\,[P]$, $(u,v)\in P$, and $(v,u)\in P$, which are clearly exclusive each other
and one of them must occur.
\end{proof}

\subsection{Height function, diameter, and the dimension of a poset}\label{ssec.poset-dimension}
This is another consequence of the previous observations.

\begin{corol}\label{cor.height-function}\em
If $n\geq 2$ then the function
$S\in {\calX}^{\preced}\setminus \{\emptyset\}\mapsto |\,\Inv(S)\,|$, extended  by a convention
$|\,\Inv(\emptyset)\,|:=-1$, is a {\em height function\/} for the lattice 
of geodetically convex subsets of\/ $\enuN$.
\end{corol}

\begin{proof}
By Theorem~\ref{thm.poset-characterization}, the discussed lattice coincides with the
lattice $({\calX}^{\preced},\subseteq)$ introduced in Section~\ref{sec.linear-extensions}.
It was shown there that it is anti-isomorphic to the poset-based lattice $({\calY}^{\preced},\subseteq)$.
Moreover, the lattice $({\calY}^{\preced},\subseteq)$ was shown in the proof of Lemma~\ref{lem.poset-lattice}(iv)
to be graded, with the height function assigning $|P\setminus\diag|$ to a poset $P$ on $N$
and ${n\choose 2}+1$ to the full relation $N\times N$. In particular, $({\calX}^{\preced},\subseteq)$
is graded, with the height function assigning the value ${n\choose 2}-|P\setminus\diag|$ to every non-empty $S={\cal L}(P)$,
which is the number of incomparable pairs in the respective poset~$P$. This value is, by Theorem~\ref{thm.poset-to-extensions}(i),
the cardinality of the set $\Inv(S)$ of inversions within $S$.
\end{proof}

Recall from Section~\ref{ssec.permut-graph} that the combinatorial equivalence of edges in the
permutohedral graph can be viewed as its particular edge-labeling, with deeper geometric meaning.
In the graphical literature, the labels on edges are typically interpreted as their ``colors".
Therefore, Corollary~\ref{cor.height-function} has the following interpretation:
the value of the height function for $({\calX}^{\preced},\subseteq)$ at non-empty
$S\in {\calX}^{\preced}$ is nothing but the number of different colors of edges within $S$.
\smallskip

The reader might have some doubts about this claim, because the elements of ${\calX}^{\preced}$ were characterized in Theorem~\ref{thm.poset-characterization} purely in graphical terms (= without using labels), while
its height function is now described in terms of particular edge labeling of the graph. Of course,
one naturally expects that the height function of the lattice $({\calX}^{\preced},\subseteq)$
should be also definable in purely graphical terms. The explanation is that the discussed {\em combinatorial
equivalence of edges\/} can be recognized/reconstructed solely on the basis of the permutohedral graph itself (as a whole).
Specifically, given an edge between $\permA\in\enuN$ and $\permB\in\enuN$, consider the set
$$
S_{\permA :\permB}:=\{\, \permC\in\enuN\,:\ \dist(\permC,\permA)<\dist(\permC,\permB)\,\}
$$
of enumerations that are closer to $\permA$ than to $\permB$. If the edge is labeled by
$\{u,v\}\subseteq N$ and $u\prec_{\permA} v$ then one can use Lemma~\ref{lem.permut-graph}(ii) to derive
that $S_{u\prec v}\subseteq S_{\permA :\permB}$. Analogously, $S_{v\prec u}\subseteq S_{\permB :\permA}$.
Since $S_{v\prec u}$ is the complement of $S_{u\prec v}$ in $\enuN$, and similarly with $S_{\permB :\permA}$ and $S_{\permA :\permB}$,
one has $S_{u\prec v}=S_{\permA :\permB}$ and $S_{v\prec u}= S_{\permB :\permA}$, meaning that
$\{S_{u\prec v},S_{v\prec u}\}=\{S_{\permA :\permB},S_{\permB :\permA}\}$ is a partition of $\enuN$, definable in purely graphical terms.
In particular, two edges have the same combinatorial label $\{u,v\}\subseteq N$ if and only if they
yield the same graphical halving of\/ $\enuN$.
\medskip

The reader may tend to think that the height function for the lattice $({\calX}^{\preced},\subseteq)$ coincides,
for $\emptyset\neq S\in {\cal X}^{\preced}$, with the
{\em diameter function\/}, which assigns the number
$$
\mbox{\rm diam}\,(S) :=\max\,\{\, \dist(\permA,\permB)\,:\ \permA,\permB\in S\,\}
\qquad \mbox{to any $\emptyset\neq S\subseteq\enuN$.}
$$
This is indeed the case if $|N|\leq 5$, but in general
one only has $\mbox{\rm diam}\,(S)\leq |\,\Inv(S)\,|$~ \cite[\S\,4]{Mas09}.

This phenomenon is closely related to the concept of the {\em dimension\/} of a poset $P$ \cite[\S\,1.1]{Wie17},
which can be equivalently defined by the formula $\mbox{\rm dim}\,(P) \,:=\, \min\,\{\, |S|\, :\ P=S^{\rGprec}\,\}$.
Note in this context that one can show that, for every pair of distinct enumerations $\permA,\permB\in\enuN$,
the set of enumerations between $\permA$ and $\permB$ belongs to ${\calX}^{\preced}$ and corresponds to a
poset of dimension 2, defined as the intersection of $T_{\permA}$ and $T_{\permB}$.
A classic result \cite{Hir51} in the theory of poset dimension says that $\mbox{\rm dim}\,(P) \leq
\lfloor \frac{n}{2}\rfloor$ for any poset $P$ on an $n$-element set $N$.
This upper bound is tight owing to a classic construction \cite{DM41} of a poset of given dimension.
Thus, the simplest example of a poset of dimension 3 exists in case $|N|=6$. It also gives
the simplest example of a poset-interpretable set $S\subseteq\enuN$ of enumerations for which
$\mbox{\rm diam}\,(S)< |\,\Inv(S)\,|$.

\begin{figure}[t]
\setlength{\unitlength}{1mm}
\begin{center}
\begin{picture}(40,24)
\put(7,5){\circle{3}}
\put(7,18){\circle{3}}
\put(20,5){\circle{3}}
\put(20,18){\circle{3}}
\put(33,5){\circle{3}}
\put(33,18){\circle{3}}
\put(7,22){\makebox(0,0){\it a}}
\put(20,22){\makebox(0,0){\it b}}
\put(33,22){\makebox(0,0){\it c}}
\put(7,1){\makebox(0,0){\it d}}
\put(20,1){\makebox(0,0){\it e}}
\put(33,1){\makebox(0,0){\it f}}
\put(8.5,16.5){\vector(1,-1){10}}
\put(31.5,16.5){\vector(-1,-1){10}}
\put(21.5,16.5){\vector(1,-1){10}}
\put(18.5,16.5){\vector(-1,-1){10}}
\put(9,16.5){\vector(2,-1){22}}
\put(31,16.5){\vector(-2,-1){22}}
\end{picture}
\end{center}
\caption{Transitive directed acyclic graph representing the poset from Example~\ref{exa.dim-poset}.\label{fig.3}}
\end{figure}

\begin{exa}\label{exa.dim-poset}\rm
Put $N:=\{a,b,c,d,e,f\}$ and $P:=\diag\cup\{\, (a,e),(a,f),(b,d),(b,f),(c,d),(c,e) \,\}$;
see Figure~\ref{fig.3} for the respective transitive directed acyclic graph.
One can show that the set $S:=P^{\lGprec}$ is the (disjoint) union of (face-associated) subsets of\/ $\enuN$:
$$
|\,abc\,|\,def\,|, \quad |\,ab\,|\,f\,|\,c\,|\,de\,|,  \quad |\,ac\,|\,e\,|\,b\,|\,d\,f\,|,
\quad |\,bc\,|\,d\,|\,a\,|\,ef\,|\,.
$$
By Theorem~\ref{thm.poset-to-extensions}(i), the inversions within $S$ are just incomparable pairs in $P$, that is, $|\,\Inv(S)\,|=9$.
On the other hand, for $\permA\in S$, if $\permA\not\in |\,ab\,|\,f\,|\,c\,|\,de\,|$ then $c\prec_{\permA} f$. This observation allows one to conclude, using Lemma~\ref{lem.permut-graph}(ii),
that any geodesic between elements of $S\setminus\, |\,ab\,|\,f\,|\,c\,|\,de\,|$ has the length at most $8$.
The implication $\permA\in S\setminus\, |\,ac\,|\,e\,|\,b\,|\,d\,f\,| ~\Rightarrow~
b\prec_{\permA} e$ and the implication $\permA\in S\setminus\, |\,bc\,|\,d\,|\,a\,|\,ef\,| ~\Rightarrow~ a\prec_{\permA} d$ lead to analogous conclusions. This, combined with the choice $\permA := |a|c|e|b|f|d|$ and $\permB := |b|c|d|a|f|e|$,
for which $\dist (\permA,\permB)=|\,\Inv [\permA,\permB ]\,|=8$
allows one to observe that $\mbox{\rm diam}\,(S)=8$.
\end{exa}

\subsection{Remarks on an alternative characterization}\label{ssec.remark-alter}
Our proof of Theorem~\ref{thm.poset-characterization} has been inspired
by the tools used in the proof of \cite[Theorem\,9]{HN13}, where a different characterization of poset-interpretable subsets
of\/ $\enuN$ was presented. It can be re-phrased as follows: a (non-empty) set $S\subseteq\enuN$
corresponds to a poset on $N$ if and only if, for every $(u,v)\in N\times N$ with $u\neq v$, exclusively one of the
following three conditions holds:
\begin{itemize}
\item[(a)] $\{u,v\}\in\Inv (S)$,
\item[(b)] $(u,v)$ belongs to the transitive closure of $\Cov (S)$,
\item[(c)] $(v,u)$ belongs to the transitive closure of $\Cov (S)$.
\end{itemize}
Observe that Corollary~\ref{cor.exclusive} claims that the above condition is necessary for $S={\cal L}(P)$,
where $P$ is a poset on $N$. The condition of this type is perhaps more suitable
from the algorithmic point of view than the condition of geodetic convexity of $S$,
when one wants to construct a poset $P$ with $\permA\in {\cal L}(P)\subseteq S$
for given $\permA\in S\subseteq\enuN$, which was the task treated in \cite[\S\,5]{HN13}.

The offered proof of the sufficiency of this condition from \cite{HN13}, however, contained three gaps.
For example, the following claim was made and used
(when Case 1.~of the proof of \cite[Thm\,9]{HN13} was discussed): if $\permA\not\in {\cal L}(P)$ for a poset~$P$ then there exists a pair $(b,a)\in P$ with $a\covA b$. This is indeed not the case: take, for instance,
$N:=\{\, a,b,c\,\}$, $P:=\diag\cup \{\,(b,a)\,\}$, and $\permA:= |a|c|b|$.

It was not clear (to the first author of the present paper) for a while whether \cite[Theorem\,9]{HN13}
is true as re-phrased above or whether there exists a counter-example of a disconnected set $S\subseteq\enuN$ satisfying exclusively (a)-(c) for some large $|N|$. Now, it looks like the statement is indeed valid as formulated above but its proof is beyond the scope of the present paper.

The reason is that proving the sufficiency requires a (correct proof of a particular)
deeper graph-theoretical result on the permutohedral graph. The needed result is, informally said,
as follows: any walk $\mathfrak{w}$ in the permutohedral graph between enumerations $\permA$ and $\permB$
can be transformed to any given geodesic $\mathfrak{g}$ between $\permA$ and $\permB$ by a series
of ``admissible" operations, which are the removals of sections of $\mathfrak{w}$ between different occurrencies
of the same node and swaps of alternative sections of $\mathfrak{w}$ of between opposite nodes on the so-called squares
and hexagons.\footnote{Squares and hexagons are particular induced subgraphs of the permutohedral graph
which correspond to two-dimensional faces of the permutohedron $\perN$ -- see Section~\ref{ssec.face-lattice}.}
Note that an algebraic analogue of such a graph-theoretical result has appeared as \cite[Thm\,10]{MPSSW09},
accompanied with arguments that are transparent only to deep specialists on Coxeter groups.
There is a plan to present a direct elementary proof of this result in the next version of \cite{Stu24}.

Here is a hint to prove the sufficiency with the help of that result. Consider an arbitrary
connectivity component $C$ of the induced subgraph $G_{S}$, where $G$ is the permutohedral graph.
The assumption on $S$ in combination with the discussed graphical result allows one to show that $C$
is geodetically convex. By Theorem~\ref{thm.poset-characterization}, a poset $P$ on $N$
exists with $C={\cal L}(P)$, and Corollary~\ref{cor.exclusive} implies that the conditions (a)-(c)
hold exclusively also for $C$. Finer considerations allow one to show that $\Inv (C)=\Inv(S)$ and $\Cov(C)=\Cov(S)$.
The application of Theorem~\ref{thm.poset-to-extensions}(iv) then gives $P=\tr(\diag\cup\Cov (C))=\tr(\diag\cup\Cov (S))$.
As the correspondence of $P\longleftrightarrow {\cal L}(P)$ is one-to-one, different components $C$
of $G_{S}$ have to coincide, implying that $S$ is connected.

\subsection{On encoding posets when a reference total order is given}\label{ssec.poset-encoding}
In the present paper, we try to propagate the opinion that when dealing with posets on a finite ground set~$N$,
one should not assume a priori that the set $N$ itself has a preferred total order. We believe that such an approach
offers simpler perspective to the audience in theoretical computer science, statistics, and combinatorics.
Nonetheless, habitual attitude of some authors dealing with this topic
is to insist on having $N$ of the form $[n]$, where $n\in {\dv N}$, that is, on having some predetermined preferred total order on~$N$. 
This approach, formalized
by the concept of a {\em labeled poset\/} from~\cite{BW91}, perhaps also has some rationale, because a fixed total order on $N$ might play the role of a kind of a coordinate system, which leads to estimates for the number of posets on~$N$.

In this section we interpret (some) results on labeled posets from \cite[\S\,6]{BW91} in our context.
It is useful to realize that, when a reference enumeration $\permD\in\enuN$ is given, which means
a total order $T_{\permD}$ on $N$ is fixed, then posets on $N$ can be identified with some intervals
in the power set of\/ $T_{\permD}\setminus\diag$, equipped with the inclusion ordering.

\begin{lemma}\label{lem.encoding-posets}\rm
Let $\permD\in\enuN$ be a fixed reference enumeration of $N$, where $|N|\geq 2$.
\begin{itemize}
\item[(i)] Given a poset $P$ on $N$, one has ${\cal L}(P)={\cal L}^{\permD}[A,B]$, where
$$
{\cal L}^{\permD}[A,B] \,:=\, \{\, \permC\in\enuN\,:\  A\subseteq T_{\permD}\setminus T_{\permC}\subseteq B\,\}
$$
with $A:=(T_{\permD}\cap P^{op})\setminus\diag$ and $B:=T_{\permD}\setminus P$, clearly
satisfying $A\subseteq B\subseteq T_{\permD}\setminus\diag$.
\item[(ii)] A set $S\subseteq\enuN$ of enumerations of $N$ has the form ${\cal L}^{\permD}[A,B]$ for
a pair $[A,B]$ of subsets of\/ $T_{\permD}\setminus\diag$ if and only if it is geodetically convex in the
permutohedral graph.
\item[(iii)] Given a non-empty geodetically convex set $S\subseteq\enuN$, the pair $[A,B]$
with $S={\cal L}^{\permD}[A,B]$ and minimal set difference $B\setminus A$ is unique.
It has the forms $[\,\,\bigcap_{\permC\in S} T_{\permD}\setminus T_{\permC}, \bigcup_{\permC\in S} T_{\permD}\setminus T_{\permC}\,]$ and $[\,(T_{\permD}\cap P^{op})\setminus\diag,T_{\permD}\setminus P\,]$, where $P$ is the unique poset on $N$ with $S={\cal L}(P)$.
\end{itemize}
\end{lemma}

The proof of this lemma, based solely on the results of the present paper, is in Appendix \ref{app.lem.encoding-posets}.
The claim (i) above means that any poset $P$ on $N$ corresponds to an interval $[A,B]$ in the poset
$(\,{\cal P}(T_{\permD}\setminus\diag)\,,\subseteq)$, where ${\cal P}(T_{\permD}\setminus\diag)$ denotes the power set of\/ $T_{\permD}\setminus\diag$. By (iii), a minimal such interval for a given $P$ is unique and
can perhaps be interpreted as a unique {\em encoding\/} of $P$, that is, yet another its cryptic description.
The claim (ii) is our reformulation of \mbox{\cite[Proposition\, 6.2]{BW91}},
attributed by Bj\"{o}rner and Wachs to Tits, this time derived without using deep group-theoretical arguments.
The claim (iii) is a kind of our re-interpretation of \cite[Proposition\, 6.3]{BW91}.
\smallskip

Surely, not every interval $[A,B]$ in $(\,{\cal P}(T_{\permD}\setminus\diag)\,,\subseteq)$ corresponds to a poset on~$N$. Nonetheless, Lemma~\ref{lem.encoding-posets} at least gives the next rough estimate for the number of posets on $N$.

\begin{corol}\label{cor.encoding-posets}\rm
The value $3^{n\choose 2}$ is an upper bound for the number of posets on $N$ with $|N|=n$.
\end{corol}

\begin{proof}
The set\, $T:=T_{\permD}\setminus\diag$ has ${n\choose 2}$ elements. Every interval $[A,B]$ in $({\cal P}(T),\subseteq)$
corresponds to a partition of\/ $T$ into three sets, namely $A$, $B\setminus A$, and $T\setminus B$.
\end{proof}

\section{Other cryptomorphic views on finite posets and conclusions}\label{sec.other-views}
Another important lattice we deal with in this paper is the lattice of {\em preposets\/} on $N$.
There are two possible interpretations of these preposets. One of them
is geometric, in terms of special cones in ${\dv R}^{N}$, as presented in \cite[\S\,3.4]{PRW08}.
The other interpretation is combinatorial, in terms of certain rings of subsets of $N$, and this
relates to well-known Birkhoff's representation theorem for finite distributive lattices \cite[\S\,III.3]{Bir95}.
Both interpretations can be introduced in terms of Galois connections, as described in
Section~\ref{ssec.app.Galois}.
Posets are special cases of preposets, for which reason the alternative views on preposets
apply to posets as well.

\subsection{Geometric view: braid cones}\label{sssec.preposet-cones}
To this end we define a binary relation $\braid$ between vectors in $X:={\dv R}^{N}$
and the elements of the Cartesian product $Y:=N\times N$. Specifically, given a vector $x\in {\dv R}^{N}$
and $(u,v)\in N\times N$, we define their incidence relation $\braid$ through the {\em comparison of the respective vector components\/}:
$$
x\,\braid\, (u,v) \quad :=\quad x_{u}\leq x_{v}\,.
$$
We will use larger black triangles $\rGbraid$ and $\lGbraid$ to
denote Galois connections based on this relation~$\braid$.
Any subset $T$ of\/ $Y=N\times N$ is assigned (by backward Galois connection) a polyhedral cone
$$
\nor_{T}:= T^{\lGbraid}= \{\, x\in {\dv R}^{N}\,:\ \forall\,(u,v)\in T~~  x_{u}\leq x_{v}\,\}\,,
$$
called the {\em braid cone\/} for $T$. In particular, any enumeration $\permA\in\enuN$ of $N$ is assigned a particular full-dimensional cone in ${\dv R}^{N}$ as the braid cone for the total order $T_{\permA}$\,:
\begin{eqnarray}
\lefteqn{\hspace*{-23mm}\nor^{\permA} \,:=~ \{\, x\in {\dv R}^{N}\,:\ x_{\permA(1)} \leq x_{\permA(2)} \leq \ldots \leq x_{\permA(n)}\,\,\}\,} \label{eq.enum-cone}\\
 &=& \{\, x\in {\dv R}^{N}\,:\ \forall\, (u,v)\in N\times N \quad u\preceq_{\permA}v ~\Rightarrow~  x_{u}\leq x_{v}\,\,\}; \nonumber
\end{eqnarray}
these cones were called {\em Weyl chambers\/} by some authors \cite{MPSSW09,PRW08}.
The next lemma claims that the method of Galois connections yields the lattice of preposets.

\begin{lemma}\label{lem.braid-lattice}\rm
Let $({\calX}^{\braid},\subseteq )$ and $({\calY}^{\braid},\subseteq )$ denote the lattices defined using the
Galois connections $\lGbraid$ and $\rGbraid$ based on the incidence relation $\braid$ introduced above.
\begin{itemize}
\item[(i)] Then $T\in {\calY}^{\braid}$ if and only if\/ $T$ is a preposet on $N$.
\item[(ii)] A preposet $T\subseteq N\times N$ is a poset if and only if its braid cone $\nor_{T}$ has
the full dimension $n$.
\item[(iii)] Given a poset $P$ on $N$, one has ${\cal L}(P)=\{\, \permA\in\enuN \,:\ \nor^{\permA}\subseteq\nor_{P}\,\}$
and $\nor_{P}=\bigcup_{\permA\in {\cal L}(P)} \nor^{\permA}$.
\end{itemize}
\end{lemma}

Our proof can be found in Appendix~\ref{app.lem.braid-lattice}.
The formulas in (iii) allow one to reconstruct ${\cal L}(P)$ on basis of $\nor_{P}$.
In fact, because Weyl chambers are isometric each other, one can perhaps interpret the ratio $\frac{|{\cal L}(P)|}{|\enuN|}$
geometrically as the relative volume of\/ $\nor_{P}\cap B$ within $B$, where $B$ is the unit ball around the origin in ${\dv R}^{N}$.

Note that the permutohedral graph can also be recognized in this geometric context. One can identify enumerations with the Weyl chambers
and the adjacency of\/ $\permA,\permB\in\enuN$ is geometrically characterized by the condition that the dimension of the intersection
$\nor^{\permA}\cap\nor^{\permB}$ is $n-1$.

\begin{remark}\label{rem.cones}\rm\,
The braid cones are also related to the order polytopes introduced by Stanley \cite{Sta86}.
Given a poset $P$ on $N$, the respective {\em order polytope\/} is defined as the intersection
$\nor_{P}\cap [0,1]^{N}$ of its braid cone with a particular $n$-dimensional cube.
Hence, given an enumeration $\permA$ of $N$, the order polytope of\/ $T_{\permA}$
is $\nor^{\permA}\cap [0,1]^{N}$.
Provided $\permA,\permB\in\enuN$ differ by an adjacent transposition \eqref{eq.adjacent-transpos},
it is an elementary observation that $\nor^{\permB}\cap [0,1]^{N}$ can be obtained from $\nor^{\permA}\cap [0,1]^{N}$ by applying the isometry of $[0,1]^{N}$ which swaps the variables $u=\permA(i)$ and $v=\permA(i+1)$.
Because these isometries preserve the Lebesgue measure, the volumes of $\nor^{\permA}\cap [0,1]^{N}$ and $\nor^{\permB}\cap [0,1]^{N}$ coincide. Since the permutohedral graph is connected, we observe by an inductive argument that the volume of $\nor^{\permA}\cap [0,1]^{N}$ does not depend on the choice of $\permA$.

Thus, the second formula in Lemma~\ref{lem.braid-lattice}(iii) allows one to conclude that the volume of the order polytope of a poset $P$ is $\frac{|{\cal L}(P)|}{|\enuN|}=\frac{|{\cal L}(P)|}{n!}$, which is an alternative derivation of \cite[Corollary\,4.2]{Sta86}.
\end{remark}

\subsection{Combinatorial view: finite topologies}\label{sssec.preposet-topologies}
Define a binary relation $\topol$ between elements of the power set $X:=\caP$ of $N$
and elements of $Y:=N\times N$. Specifically, if $D$ is a subset of $N$
and $(u,v)\in N\times N$ then we consider $D$ to be in the incidence relation $\topol$ with the pair $(u,v)$
if the set {\em $D$ respects coming $u$ before $v$\,}:
$$
D\,\topol\, (u,v) \quad := \quad [\,u\in D ~\vee~ v\not\in D\,]\,,
~~ \mbox{interpreted as $[\,v\in D \,\Rightarrow\, u\in D\,]$.}
$$
We will use larger gray triangles $\rGtopol$ and $\lGtopol$ to
denote Galois connections based on this relation~$\topol$.
Any subset $T$ of\/ $Y=N\times N$ is assigned (by backward Galois connection) a set system
$$
{\calD}_{T}:= T^{\lGtopol}= \{\, D\subseteq N\,:\ \forall\,(u,v)\in T~~  \mbox{if $v\in D$ then $u\in D$}\,\}~~\subseteq \caP\,,
$$
known as the class of {\em down sets\/} for\/ $T$, called alternatively {\em order ideals\/} of $T$ in \cite[\S\,3.1, p.\,100]{Sta97}.
In particular, if $P$ is a poset on $N$ then ${\calD}_{P}$ is nothing but what experts in (enumerative) combinatorics
denote $J(P)$ and interpret as a (finite) {\em distributive lattice\/} \cite[\S\,3.4]{Sta97}.

Note that any such a class $\calD_{T}$ satisfies $\emptyset,N\in\calD_{T}$ and is closed under set intersection and union, which means $D,E\in\calD_{T} \,\Rightarrow\, D\cap E,D\cup E\in\calD_{T}$.
Therefore, it forms a ring of sets and can be interpreted as a (finite) {\em topology on $N$}, because in case of a finite set $N$, definition of a topology on $N$ reduces to the above requirements.

Any enumeration $\permA\,=\, |\,\permA(1)|\ldots|\permA(n)\,|$ of $N$ can be represented by a
topology $\chain_{\permA}$ on $N$, the one assigned to the total order $T_{\permA}$, being the respective {\em  maximal chain\/} of subsets of~$N$:
\begin{equation}
\chain_{\permA} ~:=~ \left\{~ \bigcup_{i=1}^{j} \permA (i)\quad :\ j=0,\ldots ,n~\right\}
~~=~~ \left\{~ \emptyset ,\, \{\permA(1)\}\,,\, \{\permA(1),\permA(2)\}\,,~ \ldots ~,\,
{N} ~\right\}\,.
\label{eq.max-chain}
\end{equation}

\begin{lemma}\label{lem.topol-lattice}\rm
Let $({\calX}^{\topol},\subseteq )$ and $({\calY}^{\topol},\subseteq )$ denote the lattices defined using the
Galois connections $\lGtopol$ and $\rGtopol$ based on the incidence relation $\topol$ defined above.
\begin{itemize}
\item[(i)] Then $T\in {\calY}^{\topol}$ if and only if\/ $T$ is a preposet on $N$.
\item[(ii)] A set system $\calD\subseteq\caP$ belongs to the lattice ${\calX}^{\topol}$ if and only if
both $\emptyset,N\in\calD$ and the implication $D,E\in\calD \,\Rightarrow\, D\cap E,D\cup E\in\calD$
holds (= $\calD$ is a topology on $N$).
\item[(iii)] A preposet $T\subseteq N\times N$ is a poset if and only if
the topology $\calD_{T}$ distinguishes points: for any distinct $u,v\in N$
there is $D\in\calD_{T}$ such that either $[\,u\in D ~\&~ v\not\in D \,]$ or
$[\,v\in D ~\&~ u\not\in D \,]$. This happens exactly when there is $\permA\in\enuN$ such that
$\chain_{\permA}\subseteq\calD_{T}$.
\item[(iv)] Given a braid cone $\nor\subseteq {\dv R}^{N}$, the set system
$\calD_{\nor} :=\{\, D\subseteq N\, :\ -\chi_{D}\in\nor\,\}$ is the corresponding topology.
Given a topology $\calD$, the cone
$\nor_{\calD} := \cone(\{\chi_{N}\}\cup\{-\chi_{D}:D\in\calD\})$,
where $\cone(\nor)$ denotes the conic hull of $\nor\subseteq {\dv R}^{N}$, is the corresponding braid cone.
\item[(v)] Given a poset $P$ on $N$, one has ${\cal L}(P)=\{\, \permA\in\enuN \,:\ \chain_{\permA}\subseteq\calD_{P}\,\}$
and $\calD_{P}=\bigcup_{\permA\in {\cal L}(P)} \chain_{\permA}$.
\end{itemize}
\end{lemma}

Our proof can be found in Appendix~\ref{app.lem.topol-lattice}.
The claim (iv) describes the transfer between the geometric representation of
a preposet $T$ on~$N$, which is the braid cone $\nor_{T}$, and the respective topology $\calD_{T}$ on $N$.
Analogously, the claim (v) describes the transfer between the graphical representation ${\cal L}(P)$
of a poset $P$ on $N$ and the respective finite topology $\calD_{P}$.

Note that, by the fundamental representation theorem for finite distributive lattices, rings of subsets of a finite set (ordered by inclusion) are universal examples of such lattices. What plays substantial role in the proof of this result is a particular correspondence between finite distributive lattices and finite posets, mediated by the concept of a down set; see \cite[Corollary 1 in \S\,III.3]{Bir95}
or \cite[Theorem\,3.4.1 in \S\,3.4]{Sta97}.
This is the correspondence we discuss here, but with one important technical difference. While
distinguishing points is immaterial for the purpose of representing a (finite distributive) lattice,
it matters for the purpose of representing preposets (on the \underline{same} finite ground set).

\begin{remark}\label{rem.topologies}\rm
Let us focus on the assignment $P\mapsto \calD_{P}$ confined to posets $P$ on $N$.
A benefit of this particular transformation is that the original poset $P$ remains to be internally represented
within the distributive lattice $(\calD_{P},\subseteq)$. The point is that the variables in $N$ are in one-to-one
correspondence with the {\em join-irreducible\/} elements of this lattice, which are those sets in $\calD_{P}$
which are not the unions of their own proper subsets in $\calD_{P}$. Additionally, the comparability of
variables from $N$ in $P$ corresponds to the inclusion of the respective joint-irreducible sets in $\calD_{P}$.
Specifically, any $v\in N$ is assigned its principal order ideal $I(v):=\{\, w\in N\,:\ (w,v)\in P\,\}$ and one has $(u,v)\in P$ if and only if $I(u)\subseteq I(v)$.\footnote{To avoid misunderstanding recall that the transform $P\mapsto \calD_{P}$ itself is order-reversing when considered on the class of all posets $P$ on $N$. Here we, however, have the situation when a poset $P$ is \underline{fixed} and claim that the embedding  $v\mapsto I(v)$ of $N$ into $\calD_{P}$ is order-preserving.}

The transformation $P\mapsto \calD_{P}$ of a poset $P$ on $N$ is also one of standard constructions in combinatorics
because it leads to more elaborate methods to count linear extensions of a poset.
Indeed, Lemma~\ref{lem.topol-lattice}(v) allows one to identify {\em linear extensions\/} of\/ $P$
on basis of\/ $\calD_{P}$ and, thus, this two-step procedure $P\mapsto \calD_{P}\mapsto {\cal L}(P)$ forms the basis for common algorithms to compute the number of linear extensions of\/ $P$, as indicated in \cite[\S\,3.5]{Sta97}. A particular way of representing the lattice $(\calD_{P},\subseteq)$ in the form of a grid in ${\dv N}^{k}$ for some $k$ is described there and the linear extensions of $P$ are identified with certain ``increasing" paths in this grid. The problem of counting linear extensions of $P$ is thus transformed to the problem of counting certain types of permutations.

According to the explanation in \cite[p.\,110]{Sta97}, it follows from a classic result \cite{Dil50} by Dilworth
that the (poset) dimension of $(\calD_{P},\subseteq)$ is the {\em width\/} of the poset $P$, which is the maximal
cardinality of an anti-chain in $P$, being the least $k$ such that $N$ can be written as the union of\/ $k$ chains
in~$P$. Thus, the width of $P$ seems to be the crucial complexity measure for computing $|{\cal L}(P)|$.
\end{remark}

\subsection{Conclusions}\label{ssec.conclusions}
In the present paper, we have described some cryptomorphic views on posets on a non-empty finite set of variables.
A contribution of the second author is a web page
\begin{quote}
{\tt  http://gogo.utia.cas.cz/posets/}
\end{quote}
which allows the user to pass through different combinatorial representatives of posets
on a variable set $N$ of small cardinality. Specifically, besides the basic representation of a poset $P$ in the form
of a binary relation on $N$, it can be represented by the respective transitive directed acyclic graph over $N$, by its
directed Hasse diagram, by the set of its linear extensions ${\cal L}(P)$, and by the respective
finite topology $\calD_{P}$, that is, in the form of a finite distributive lattice.

We believe that the cryptomorphic views on finite posets will appear to be useful in context
of combinatorial description of faces of the supermodular cone \cite{Stu24} by means of collections of finite posets,
called the {\em fans of posets} \cite[\S\,3.3]{PRW08}.
Note that the fans of posets have already appeared in the literature on (learning) graphical statistical models, where they were related to polytopes known as {\em generalized permutohedra}.
Specifically, in \cite[\S\,8]{MUWY18}, an example was given in which a Bayesian network over a 4-element ground set was represented both as a collection of posets and as a collection of corresponding sets of enumerations.

We also hope that our main result on graphical characterization of sets of linear extensions of posets
as geodetically convex sets might find future applications in alternative algorithms to deal with the
poset cover problem, since it is conceptually simpler than \cite[Theorem\,9]{HN13}.


\subsubsection*{Addendum (from the first author)}
I can hardly claim that the theme of posets was in the center of research
activity of Henry Wynn. Nonetheless, he was open and willing to discuss
this topic and its possible connections to Gr\"{o}bner bases applications
and we did so together on several occasions of conferences on algebraic statistics.
I am also indebted to Victor Reiner for drawing my attention to the results
presented in \cite{BW91,Tit74}, which are highly relevant to the topic of the present paper.
And, finally, both authors are indebted to the reviewers and the editors for their positive attitude and
their recommendations, which improved the quality and readability of the paper.

\appendix

\section{Proof of Lemma~\ref{lem.poset-lattice}}\label{app.lem.poset-lattice}
For the reader's convenience we restate the lemmas and theorems to be proved throughout the appendix.
\medskip

\noindent{\bfseries Lemma~\ref{lem.poset-lattice}:} ~~\rm
Let $({\calX}^{\preced},\subseteq )$ and $({\calY}^{\preced},\subseteq )$ denote the finite lattices defined using the Galois connections $\lGprec$ and $\rGprec$ based on the incidence relation $\preced$ between $X$ and $Y$ introduced in Section~\ref{sec.linear-extensions}.
\begin{itemize}
\item[(i)] If $n=|N|\geq 2$ then\, [\,$T\in {\calY}^{\preced}$ and\, $T\subset N\times N$\,]
if and only if\/ $T$ is a poset on $N$.
\item[(ii)] One has [\,$S\in {\calX}^{\preced}$ and $S\neq\emptyset$\,] if and only if
a (uniquely determined) poset $(N,\preceq)$ exists such that $S=\{\,\permA\in\enuN\,:\
\permA_{-1}(u)<\permA_{-1}(v)~~ \mbox{whenever $u\prec v$}\,\}$, that is,  $S={\cal L}(P)$ for a poset~$P$ on $N$.
\item[(iii)] Any $S\subseteq\enuN$ representing a non-empty face of\/ $\Pi (N)$
(as in Section\,\ref{ssec.face-lattice})  belongs to ${\calX}^{\preced}$.
\item[(iv)] The lattice $({\calX}^{\preced},\subseteq )$ is
atomistic, coatomistic, and graded.
\end{itemize}
Coatoms in ${\calX}^{\preced}$ are precisely the sets $S_{u\prec v}\,:=\,\{\,\permA\in\enuN \,:\ u\prec_{\permA}v\,\}$,
where $u,v\in N$, $u\neq v$.

\begin{proof}
For the necessity in (i) assume $T\in{\calY}^{\preced}$ and\/ $T\neq N\times N$. Therefore,
$$
T=S^{\rGprec}:=\{\, (u,v)\in N\times N :\
\forall\,\permA\in S~~ \permA_{-1}(u)\leq\permA_{-1}(v)\,\}
$$
for some $S\subseteq\enuN$ (see Section~\ref{ssec.app.Galois}).
The reflexivity and transitivity of\/ $T$ follows directly from its definition on the basis of $S$.
Then $T$ is anti-symmetric as otherwise $S=\emptyset$ gives $T=N\times N$.

For the sufficiency in (i) assume that $T\subset N\times N$ is a poset. By Fact~\ref{fct.posets}
in Section \ref{ssec.trans-closure}, the relation
$T\setminus\diag$ is then acyclic and coincides with the set of arrows of
a transitive directed acyclic graph $G$ over $N$ (see Section~\ref{ssec.poset-as-DAG}). Hence, the set
\begin{eqnarray*}
S\,:=\, T^{\lGprec} &=& \{\, \permA\in\enuN\,:\ \forall\,(u,v)\in T~~ \permA_{-1}(u)\leq\permA_{-1}(v)\,\}\\
&=& \{\, \permA\in\enuN\,:\  \permA_{-1}(a)<\permA_{-1}(b)\quad
\mbox{whenever $a\to b$ in $G$}\,\}
\end{eqnarray*}
of enumerations of $N$ which are in concordance with $G$ is non-empty. Evidently, $T\subseteq S^{\rGprec}$ and to show $S^{\rGprec}\subseteq T$ one needs, for any $(u,v)\not\in T$, to find an enumeration $\permA\in S$ with $\permA_{-1}(v)<\permA_{-1}(u)$.
Since $S\neq\emptyset$, this is evident in case $v\to u$ in $G$.
In the case of a pair $(u,v)$ of non-adjacent distinct nodes of $G$ one can consider
the set $A\subseteq N$ of $a\in N$ with $(a,v)\in T$. Then $v\in A$, $u\not\in A$, and
the induced subgraph $G_{A}$ is acyclic. Hence, an enumeration of $A$ exists which is
in concordance with $G_{A}$ and this can be extended to an enumeration $\permA$ of $N$ in concordance with~$G$.
Thus, $\permA\in S$ with $\permA_{-1}(v)<\permA_{-1}(u)$ was obtained.
Altogether, we have shown $T=S^{\rGprec}$, which means that $T\in {\calY}^{\preced}$.

To prove (ii), note that it holds trivially in case $n=1$. In case $n\geq 2$ apply (i) and realize
that the correspondence $T\in {\calY}^{\preced}\mapsto T^{\lGprec}={\cal L}(T)\in {\calX}^{\preced}$
is an anti-isomorphism (see Section~\ref{ssec.app.Galois}) and maps bijectively the partial orders on $N$ onto non-empty sets in ${\calX}^{\preced}$.

Note that (iii) holds if $n=1$. If $n\geq 2$ then take a non-empty face of $\perN$ represented by an ordered partition
$|\,A_{1}|\ldots|A_{m}\,|$, $m\geq 1$, of $N$ (see Section~\ref{ssec.face-lattice}). It defines a partial order
$$
T:=\diag\cup \{\, (u,v)\in N\times N\,:\ u\in A_{i} ~\mbox{and}~ v\in A_{j}~~ \mbox{for some $i<j$}\,\,\}
$$
on $N$. Easily, ${\cal L}(T)$ coincides with the set of enumerations
consistent with $|\,A_{1}|\ldots|A_{m}\,|$.
\smallskip

To prove (iv), note that it holds trivially in case $n=1$. In case $n\geq 2$ consider (any) $\permA\in\enuN$ and observe that $(T_{\permA})^{\lGprec}=\{\,\permA\,\}$.
Hence, the atoms of $({\calX}^{\preced},\subseteq )$ are the singleton subsets of\/ $\enuN$. Since any $S\in {\calX}^{\preced}$ is the union of singletons, the lattice $({\calX}^{\preced},\subseteq )$ is atomistic.

The fact that $({\calX}^{\preced},\subseteq )$ is coatomistic follows from the fact that $({\calY}^{\preced},\subseteq )$ is atomistic. To show the latter fact realize that, for
every $(u,v)\in (N\times N)\setminus\diag$, the set $\diag\cup \{\,(u,v)\,\}$ is a poset on~$N$.
Thus, using (i), the atoms of $({\calY}^{\preced},\subseteq )$ correspond to
singleton subsets of $(N\times N)\setminus\diag$. This allows one to observe, again
using (i), that $({\calY}^{\preced},\subseteq )$ is atomistic.

Note that $({\calX}^{\preced},\subseteq )$ is graded if and only if
$({\calY}^{\preced},\subseteq )$ is graded. To verify the latter fact one can use the characterization of\/ ${\calY}^{\preced}$ from (i).
The anti-isomorphism of $({\calX}^{\preced},\subseteq )$ and
$({\calY}^{\preced},\subseteq )$ implies that the coatoms of $({\calY}^{\preced},\subseteq )$
are the total orders $T_{\permA}$, $\permA\in\enuN$, where every
$T_{\permA}\setminus\diag$ has the cardinality ${n\choose 2}$.
As explained in Section~\ref{ssec.covering-relation},
to reach our goal it is enough to show that any maximal
chain in $({\calY}^{\preced}\setminus\{N\times N\},\subseteq )$
has the length ${n\choose 2}$, that is, it has ${n\choose 2}+1$ elements.

This follows by repeated application of the following {\em sandwich principle}. Given two posets
$P^{\prime}\subset P$ on $N$ there exists a poset $P^{\prime\prime}$ on $N$ with
$P^{\prime}\subseteq P^{\prime\prime}\subset P$ and $|P\setminus P^{\prime\prime}|=1$.
To verify the principle realize that there exists $(u,v)\in N\times N$ with $u\cover v\,\,[P]$
and $(u,v)\not\in P^{\prime}$.\\[0.4ex]
Indeed, otherwise the covering relation $\cover$ for $P$ is contained in $P^{\prime}$ and its transitive closure, which is $P\setminus\diag$
(by Fact~\ref{fct.covering} in Section~\ref{ssec.poset-as-DAG}), is also contained, yielding a contradictory fact $P\subseteq P^{\prime}$.\\[0.4ex]
Take a pair $(u,v)$ with $u\cover v\,\,[P]$ and $(u,v)\not\in P^{\prime}$ and put $P^{\prime\prime}:=P\setminus\{\,(u,v)\,\}$.
The transitive closure of\/ $P^{\prime\prime}$ is contained in (a transitive relation) $P$. Nevertheless, $(u,v)$ cannot be in the transitive closure of\/ $P^{\prime\prime}$, because this contradicts the assumption that $u\cover v\,\,[P]$. Thus, $P^{\prime\prime}\setminus\diag$ is both transitive and irreflexive, which verifies the sandwich principle (see Section~\ref{ssec.trans-closure}).

The claim about coatoms in ${\calX}^{\preced}$ follows from the description of atoms in $({\calY}^{\preced},\subseteq)$:
these have the form $A_{(u,v)}=\diag\cup \{\, (u,v)\,\}$ for $u,v\in N$, $u\neq v$, and $A_{(u,v)}^{\lGprec}=\{\, (u,v)\,\}^{\lGprec}=S_{u\prec v}$.
\end{proof}

\section{Proof of Lemma~\ref{lem.permut-graph}}\label{app.lem.permut-graph}

\noindent{\bfseries Lemma~\ref{lem.permut-graph}:} ~~\rm
Consider $\permA,\permB\in\enuN$. Then
\begin{itemize}
\item[(i)] $\dist(\permA,\permB)=|\,\Inv[\permA,\permB]\,|$ is the number of inversions
between $\permA$ and $\permB$. In particular, $\permA$ and~$\permB$ are adjacent
in the permutohedral graph if and only if they differ by one inversion only. The label on the edge
between $\permA$ and $\permB$ is then the only inversion between them.
\item[(ii)]
A walk in the permutohedral graph is a geodesic (between $\permA$ and $\permB$) if and only if no label 
on its edges is repeated. If this is the case then the set of labels on the geodesic coincides
with the set $\Inv[\permA,\permB]$ of inversions between $\permA$ and $\permB$.
\item[(iii)]
An enumeration $\permC\in\enuN$ is between $\permA$ and $\permB$ if and only if $\Inv[\permA,\permC]\cap\Inv[\permC,\permB]=\emptyset$.
Another equivalent condition is $T_{\permA}\cap T_{\permB}\subseteq T_{\permC}$, that is,
$[\, u\prec_{\permA} v ~\&~  u\prec_{\permB} v\,] ~\Rightarrow ~ u\prec_{\permC} v$.
\end{itemize}

\begin{proof}
Assume $n\geq 2$ as otherwise all claims are trivial.
Let us start with the (particular) observation in (i) that $\permA,\permB\in\enuN$ are adjacent
if and only if they differ by one inversion only. The necessity follows immediately from \eqref{eq.adjacent-transpos}
(Section~\ref{ssec.permut-graph}).
For the sufficiency assume that $\{u,v\}\in\Inv[\permA,\permB]$ is unique, meaning
that any other pair of elements of $N$ has the same mutual order in $\permA$ and in~$\permB$.
Assume specifically that $u\prec_{\permA}v$ and $v\prec_{\permB}u$ as otherwise we exchange
$\permA$ for $\permB$.
Thus, one has $u=\permA(i)$ and $v=\permA(j)$ with $i<j$.
One has $j=i+1$ as otherwise there is $t\in N$ with $u\prec_{\permA}t\prec_{\permA}v$ implying
a contradictory conclusion $u\prec_{\permB}t\prec_{\permB}v\prec_{\permB}u$.
Moreover, for any $t\in N\setminus\{u,v\}$, one has $t\prec_{\permA}u \,\Leftrightarrow\, t\prec_{\permB}u$, and
$v\prec_{\permA}t \,\Leftrightarrow\, v\prec_{\permB}t$. These facts enforce $v=\permB(i)$ and $u=\permB(i+1)$.
The fact that the mutual orders in $\permA$ and $\permB$ coincide among $t\in N$ with $t\prec_{\permA}u$
and also among $t\in N$ with $v\prec_{\permA}t$ then implies $\permA(k)=\permB(k)$ for $k\in[n]\setminus\{i,i+1\}$.
Thus, \eqref{eq.adjacent-transpos} holds.
Note that the label $\{u,v\}$ on the edge between enumerations $\permA$ and $\permB$ is just the only
inversion between them.

The latter observation implies that, for every walk between enumerations $\permA,\permB\in\enuN$
and every $\{u,v\}\in\Inv[\permA,\permB]$, there is an edge of the walk labeled by $\{u,v\}$, as otherwise
$u$ and $v$ would have the same mutual order in $\permA$ and in $\permB$. Hence, $|\,\Inv[\permA,\permB]\,|$
is a lower bound for the length of a walk between $\permA$ and $\permB$.

Therefore, to verify (i) it is enough to prove, for $\permA,\permB\in\enuN$, by induction on $|\,\Inv[\permA,\permB]\,|$, that there exists a walk between them which has the length $|\,\Inv[\permA,\permB]\,|$.
If there is no inversion between $\permA$ and $\permB$ then
$\permA(1)\prec_{\permA} \permA(2)\, \ldots \prec_{\permA} \permA(n)$ gives
$\permA(1)\prec_{\permB} \permA(2)\, \ldots \prec_{\permB} \permA(n)$ and $\permA=\permB$,
that is, a walk of the length $0=|\,\Inv[\permA,\permB]\,|$ exists.
If $|\,\Inv[\permA,\permB]\,|\geq 1$ then necessarily $\permA\neq\permB$ and
$\exists\, 1\leq i<n$ such that one has $\permB(i+1)\prec_{\permA} \permB(i)$, for
otherwise $\permB(1)\prec_{\permA} \permB(2)\, \ldots \prec_{\permA} \permB(n)$ would give
a contradictory conclusion $\permB=\permA$. Define $\permC\in\enuN$ as the enumeration
obtained from $\permB$ by the respective adjacent transposition: $\permC(i)=\permB(i+1)=:u$,
$\permC(i+1)=\permB(i)=:v$, and $\permC(k)=\permB(k)$ for $k\in[n]\setminus\{i,i+1\}$.
Since $\{u,v\}$ is an inversion between $\permA$ and $\permB$,
one easily observes $\Inv[\permA,\permC]=\Inv[\permA,\permB]\setminus\{\, \{u,v\}\,\}$.
The induction premise implies the existence of a walk between $\permA$ and $\permC$ of the
length $|\,\Inv[\permA,\permC]\,|$ which can be prolonged by the edge between $\permC$ and $\permB$.
That verifies the induction step.

To prove the necessity in (ii) consider a geodesic $\mathfrak{g}$ between $\permA,\permB\in\enuN$
and assume for a contradiction that a label $\{u,v\}\subseteq N$ is repeated on its edges.
The last claim in (i) implies that the label cannot be repeated on consecutive edges of $\mathfrak{g}$:
if an edge between $\permD\in\enuN$ and $\permE\in\enuN$ has the same label as the
consecutive edge between $\permE$ and $\permC\in\enuN$ then $\permD=\permC$ and $\mathfrak{g}$ can be shortened,
contradicting the assumption that it is a geodesic.

Therefore, one can assume without loss of generality
that the first occurrence $\{u,v\}$ on the way from $\permA$ to $\permB$ is
between enumerations $\permD$ and $\permE$, with $\permD$ closer to $\permA$, and the next
occurrence of $\{u,v\}$ is between enumerations $\permC$ and $\permF$, with $\permF$ closer to $\permB$
(see Figure~\ref{fig.1}(a) for illustration).

\begin{figure}[t]
\setlength{\unitlength}{0.8mm}
\begin{center}
\begin{picture}(79,33)
\put(5,20){\circle{3}}
\put(20,20){\circle{3}}
\put(20,7){\color{lgrey}\circle*{3}}
\put(33,7){\color{lgrey}\circle*{3}}
\put(46,7){\color{lgrey}\circle*{3}}
\put(59,7){\color{lgrey}\circle*{3}}
\put(59,20){\circle{3}}
\put(74,13.5){\circle{3}}
\put(5,24){\makebox(0,0){\footnotesize $\permA$}}
\put(20,3){\makebox(0,0){\footnotesize $\permE\in S_{v\prec u}$}}
\put(39.5,3){\makebox(0,0){\small $\mathfrak{g}^{\prime}$}}
\put(59,3){\makebox(0,0){\footnotesize $\permC\in S_{v\prec u}$}}
\put(20,24){\makebox(0,0){\footnotesize $\permD\in S_{u\prec v}$}}
\put(59,24){\makebox(0,0){\footnotesize $\permF\in S_{u\prec v}$}}
\put(74,17.5){\makebox(0,0){\footnotesize $\permB$}}
\put(39.5,31){\makebox(0,0){\rm (a)}}
\put(2,13.5){\line(1,0){2}}
\put(5,13.5){\line(1,0){2}}
\put(8,13.5){\line(1,0){2}}
\put(11,13.5){\line(1,0){2}}
\put(14,13.5){\line(1,0){2}}
\put(17,13.5){\line(1,0){2}}
\put(20,13.5){\line(1,0){2}}
\put(23,13.5){\line(1,0){2}}
\put(26,13.5){\line(1,0){2}}
\put(29,13.5){\line(1,0){2}}
\put(32,13.5){\line(1,0){2}}
\put(35,13.5){\line(1,0){2}}
\put(38,13.5){\line(1,0){2}}
\put(41,13.5){\line(1,0){2}}
\put(44,13.5){\line(1,0){2}}
\put(47,13.5){\line(1,0){2}}
\put(50,13.5){\line(1,0){2}}
\put(53,13.5){\line(1,0){2}}
\put(56,13.5){\line(1,0){2}}
\put(59,13.5){\line(1,0){2}}
\put(62,13.5){\line(1,0){2}}
\put(20,18.1){\line(0,-1){9.2}}
\put(59,18.1){\line(0,-1){9.2}}
\put(8.5,20){\circle*{0.4}}
\put(10.5,20){\circle*{0.4}}
\put(12.5,20){\circle*{0.4}}
\put(14.5,20){\circle*{0.4}}
\put(16.5,20){\circle*{0.4}}
\put(62.5,18.75){\circle*{0.4}}
\put(64.5,17.75){\circle*{0.4}}
\put(66.5,16.75){\circle*{0.4}}
\put(68.5,15.75){\circle*{0.4}}
\put(70.5,14.75){\circle*{0.4}}
\put(36.5,7){\color{lgrey}\circle*{0.8}}
\put(39.5,7){\color{lgrey}\circle*{0.8}}
\put(42.5,7){\color{lgrey}\circle*{0.8}}
\thicklines
\put(21.9,7){\color{lgrey}\line(1,0){9.2}}
\put(47.9,7){\color{lgrey}\line(1,0){9.2}}
\end{picture}
~~~~~
\begin{picture}(79,33)
\put(5,20){\circle{3}}
\put(20,20){\color{lgrey}\circle*{3}}
\put(33,20){\color{lgrey}\circle*{3}}
\put(46,20){\color{lgrey}\circle*{3}}
\put(59,20){\color{lgrey}\circle*{3}}
\put(20,7){\circle{3}}
\put(33,7){\circle{3}}
\put(46,7){\circle{3}}
\put(59,7){\circle{3}}
\put(74,13.5){\circle{3}}
\put(5,24){\makebox(0,0){\footnotesize $\permA$}}
\put(20,3){\makebox(0,0){\footnotesize $\permE$}}
\put(39.5,24){\makebox(0,0){\small $\mathfrak{g}^{\prime\prime}$}}
\put(59,3){\makebox(0,0){\footnotesize $\permC$}}
\put(20,24){\makebox(0,0){\footnotesize $\permD$}}
\put(59,24){\makebox(0,0){\footnotesize $\permF$}}
\put(74,17.5){\makebox(0,0){\footnotesize $\permB$}}
\put(44,13.5){\makebox(0,0){\footnotesize $\trans_{uv}$}}
\put(39.5,31){\makebox(0,0){\rm (b)}}
\put(2,13.5){\line(1,0){2}}
\put(5,13.5){\line(1,0){2}}
\put(8,13.5){\line(1,0){2}}
\put(11,13.5){\line(1,0){2}}
\put(14,13.5){\line(1,0){2}}
\put(17,13.5){\line(1,0){2}}
\put(20,13.5){\line(1,0){2}}
\put(23,13.5){\line(1,0){2}}
\put(26,13.5){\line(1,0){2}}
\put(50,13.5){\line(1,0){2}}
\put(53,13.5){\line(1,0){2}}
\put(56,13.5){\line(1,0){2}}
\put(59,13.5){\line(1,0){2}}
\put(62,13.5){\line(1,0){2}}
\put(20,18.1){\line(0,-1){9.2}}
\put(59,18.1){\line(0,-1){9.2}}
\put(39.5,16.5){\vector(0,-1){6}}
\put(39.5,10.5){\vector(0,1){6}}
\put(8.5,20){\circle*{0.4}}
\put(10.5,20){\circle*{0.4}}
\put(12.5,20){\circle*{0.4}}
\put(14.5,20){\circle*{0.4}}
\put(16.5,20){\circle*{0.4}}
\put(62.5,18.75){\circle*{0.4}}
\put(64.5,17.75){\circle*{0.4}}
\put(66.5,16.75){\circle*{0.4}}
\put(68.5,15.75){\circle*{0.4}}
\put(70.5,14.75){\circle*{0.4}}
\put(36.5,7){\circle*{0.6}}
\put(39.5,7){\circle*{0.6}}
\put(42.5,7){\circle*{0.6}}
\put(36.5,20){\color{lgrey}\circle*{0.8}}
\put(39.5,20){\color{lgrey}\circle*{0.8}}
\put(42.5,20){\color{lgrey}\circle*{0.8}}
\put(21.9,7){\line(1,0){9.2}}
\put(47.9,7){\line(1,0){9.2}}
\thicklines
\put(21.9,20){\color{lgrey}\line(1,0){9.2}}
\put(47.9,20){\color{lgrey}\line(1,0){9.2}}
\end{picture}
\caption{A picture illustrating the proof of Lemma~\ref{lem.permut-graph}.\label{fig.1}}
\end{center}
\end{figure}

Assume specifically $u\prec_{\permA} v$, which enforces $\permD,\permF\in S_{u\prec v}$ and
$\permE,\permC\in S_{v\prec u}$. Let $\mathfrak{g}^{\prime}$ be the section of $\mathfrak{g}$
between $\permE$ and $\permC$. We know $\mathfrak{g}^{\prime}\subseteq S_{v\prec u}$ and, because the map
$\permG\mapsto\trans_{uv}\cdot\permG$ is an automorphism of the permutohedral graph (see Section~\ref{ssec.permut-graph}), the section $\mathfrak{g}^{\prime}$
has a copy $\mathfrak{g}^{\prime\prime}$ in $S_{u\prec v}$, which is between $\permD$ and $\permF$;
see Figure~\ref{fig.1}(b).
The section of $\mathfrak{g}$ between $\permD$ and $\permF$ can be replaced by $\mathfrak{g}^{\prime\prime}$,
yielding a walk of a shorter length and contradicting the assumption that $\mathfrak{g}$ is a geodesic.

To prove the sufficiency in (ii) assume that $\mathfrak{g}$ is a walk between $\permA$ and $\permB$
in which no label is repeated. We have already observed that any inversion $\{u,v\}\in\Inv[\permA,\permB]$
must occur as a label on $\mathfrak{g}$. Nevertheless, no two-element set $\{u,v\}\subseteq N$ which
is {\em not an inversion\/} between $\permA$ and $\permB$ can be a label on $\mathfrak{g}$:
indeed, otherwise, since only one edge of $\mathfrak{g}$ is labeled by $\{u,v\}$ then, the mutual orders
of $u$ and $v$ in $\permA$ and $\permB$ must differ, which contradicts the assumption
that $\{u,v\}$ is not an inversion between $\permA$ and $\permB$. Thus, we have verified
both that the length of $\mathfrak{g}$ is $|\,\Inv[\permA,\permB]\,|=\dist(\permA,\permB)$
and the claim that $\Inv[\permA,\permB]$ is the set of labels on $\mathfrak{g}$.
\smallskip

The first claim in (iii) is implied by (ii) as follows. If $\permC$ belongs to a geodesic between $\permA$
and $\permB$ then the labels on it between $\permA$ and $\permC$ must differ from the labels between $\permC$ and $\permB$.
Conversely, if $\Inv[\permA,\permC]$ and $\Inv[\permC,\permB]$ are disjoint then the concatenation
of a geodesic between $\permA$ and $\permC$ with a geodesic between $\permC$ and $\permB$ yields a geodesic
between $\permA$ and $\permB$.

Assuming $\Inv[\permA,\permC]\cap\Inv[\permC,\permB]=\emptyset$, the conditions
$u\preceq_{\permA}v$ and $u\preceq_{\permB}v$ together imply $u\preceq_{\permC}v$
as otherwise $\{u,v\}$ belongs to the intersection of the sets of inversions.
Conversely, if $T_{\permA}\cap T_{\permB}\subseteq T_{\permC}$ then
any hypothetic simultaneous inversion $\{u,v\}$ in $\Inv[\permA,\permC]\cap\Inv[\permC,\permB]$
must be compared equally in $\permA$ and $\permB$, and, therefore, in $\permC$, which is a contradiction.
\end{proof}

\section{Proof of Theorem~\ref{thm.poset-characterization}}\label{app.thm.poset-characterization}

\noindent{\bfseries Theorem~\ref{thm.poset-characterization}:} ~~\rm
If $n=|N|\geq 2$ then, given $S\subseteq\enuN$, $S\in {\calX}^{\preced}$ if and only if $S$ is geodetically convex.
In particular, $S={\cal L}(P)$ for a poset $P$ on $N$ if and only if $S$ is a non-empty geodetically convex set of enumerations of $N$.

\begin{proof}
{\bf I.}\, To verify the necessity use Lemma~\ref{lem.poset-lattice}(iv) saying that
$({\calX}^{\preced},\subseteq)$ is a coatomistic lattice. Since the intersection of geodetically convex
sets is geodetically convex and so is $S=\enuN$ it is enough to verify that
the coatoms in $({\calX}^{\preced},\subseteq)$ are geodetically convex. These are sets of the form
$S_{u\prec v}$ for pairs $(u,v)\in N\times N$, $u\neq v$ (see Lemma~\ref{lem.poset-lattice}). We apply Lemma~\ref{lem.permut-graph}(ii):
given a geodesic between $\permA,\permB\in S_{u\prec v}$, the fact $\{u,v\}\not\in \Inv[\permA,\permB]$
implies that $\{u,v\}$ is not a label on any edge of the geodesic, and, thus, every node on the geodesic
belongs to $S_{u\prec v}$.

\begin{figure}[t]
\setlength{\unitlength}{0.8mm}
\begin{center}
\begin{picture}(79,33)
\put(7,7){\circle{3}}
\put(7,20){\color{lgrey}\circle*{3}}
\put(20,20){\color{lgrey}\circle*{3}}
\put(33,20){\color{lgrey}\circle*{3}}
\put(46,20){\color{lgrey}\circle*{3}}
\put(59,7){\circle{3}}
\put(74,7){\circle{3}}
\put(7,3){\makebox(0,0){\footnotesize $\permB\not\in S$}}
\put(7,24){\makebox(0,0){\footnotesize $\permA\in S$}}
\put(26.5,24){\makebox(0,0){\small $\mathfrak{g}^{\prime}$}}
\put(46,24){\makebox(0,0){\footnotesize $\permD\in S$}}
\put(46,3){\makebox(0,0){\footnotesize $\permE\in S$}}
\put(74,3){\makebox(0,0){\footnotesize $\permC\in S$}}
\put(4,13.5){\makebox(0,0){\scriptsize $uv$}}
\put(49,13.5){\makebox(0,0){\scriptsize $uv$}}
\put(67,23){\makebox(0,0){\small $S_{u\prec v}$}}
\put(38,31){\makebox(0,0){\rm (a)}}
\put(62.5,7){\circle*{0.4}}
\put(64.5,7){\circle*{0.4}}
\put(66.5,7){\circle*{0.4}}
\put(68.5,7){\circle*{0.4}}
\put(70.5,7){\circle*{0.4}}
\put(11,13.5){\line(1,0){2}}
\put(14,13.5){\line(1,0){2}}
\put(17,13.5){\line(1,0){2}}
\put(20,13.5){\line(1,0){2}}
\put(23,13.5){\line(1,0){2}}
\put(26,13.5){\line(1,0){2}}
\put(29,13.5){\line(1,0){2}}
\put(32,13.5){\line(1,0){2}}
\put(35,13.5){\line(1,0){2}}
\put(38,13.5){\line(1,0){2}}
\put(41,13.5){\line(1,0){2}}
\put(53,13.5){\line(1,0){2}}
\put(56,13.5){\line(1,0){2}}
\put(59,13.5){\line(1,0){2}}
\put(62,13.5){\line(1,0){2}}
\put(65,13.5){\line(1,0){2}}
\put(68,13.5){\line(1,0){2}}
\put(71,13.5){\line(1,0){2}}
\put(7,18.1){\line(0,-1){9.2}}
\put(47.9,7){\line(1,0){9.2}}
%
\put(23.5,20){\color{lgrey}\circle*{0.8}}
\put(26.5,20){\color{lgrey}\circle*{0.8}}
\put(29.5,20){\color{lgrey}\circle*{0.8}}
\thicklines
\put(8.9,20){\color{lgrey}\line(1,0){9.2}}
\put(34.9,20){\color{lgrey}\line(1,0){9.2}}
\put(46,18.1){\color{hgrey}\line(0,-1){9.2}}
\put(7,20){\circle{3}}
\put(46,7){\circle{3}}
\end{picture}
~~~~~
\begin{picture}(79,33)
\put(20,20){\circle{3}}
\put(33,20){\circle{3}}
\put(46,20){\circle{3}}
\put(7,7){\color{lgrey}\circle*{3}}
\put(20,7){\color{lgrey}\circle*{3}}
\put(33,7){\color{lgrey}\circle*{3}}
\put(46,7){\color{lgrey}\circle*{3}}
\put(59,7){\circle{3}}
\put(74,7){\circle{3}}
\put(7,3){\makebox(0,0){\footnotesize $\permB$}}
\put(7,24){\makebox(0,0){\footnotesize $\permA\in S$}}
\put(26.5,3){\makebox(0,0){\small $\mathfrak{g}^{\prime\prime}$}}
\put(46,24){\makebox(0,0){\footnotesize $\permD\in S$}}
\put(46,3){\makebox(0,0){\footnotesize $\permE$}}
\put(74,3){\makebox(0,0){\footnotesize $\permC\in S$}}
\put(31,13.5){\makebox(0,0){\footnotesize $\trans_{uv}$}}
\put(67,23){\makebox(0,0){\small $S_{u\prec v}$}}
\put(38,31){\makebox(0,0){\rm (b)}}
\put(62.5,7){\circle*{0.4}}
\put(64.5,7){\circle*{0.4}}
\put(66.5,7){\circle*{0.4}}
\put(68.5,7){\circle*{0.4}}
\put(70.5,7){\circle*{0.4}}
%
\put(11,13.5){\line(1,0){2}}
\put(14,13.5){\line(1,0){2}}
\put(17,13.5){\line(1,0){2}}
\put(20,13.5){\line(1,0){2}}
\put(38,13.5){\line(1,0){2}}
\put(41,13.5){\line(1,0){2}}
\put(50,13.5){\line(1,0){2}}
\put(53,13.5){\line(1,0){2}}
\put(56,13.5){\line(1,0){2}}
\put(59,13.5){\line(1,0){2}}
\put(62,13.5){\line(1,0){2}}
\put(65,13.5){\line(1,0){2}}
\put(68,13.5){\line(1,0){2}}
\put(71,13.5){\line(1,0){2}}
%
\put(46,18.1){\line(0,-1){9.2}}
\put(47.9,7){\line(1,0){9.2}} 
\put(8.9,20){\line(1,0){9.2}}
\put(34.9,20){\line(1,0){9.2}}
\put(26.5,16.5){\vector(0,-1){6}}
\put(26.5,10.5){\vector(0,1){6}}
\put(23.5,20){\circle*{0.6}}
\put(26.5,20){\circle*{0.6}}
\put(29.5,20){\circle*{0.6}}
\put(23.5,7){\color{lgrey}\circle*{0.8}}
\put(26.5,7){\color{lgrey}\circle*{0.8}}
\put(29.5,7){\color{lgrey}\circle*{0.8}}
\thicklines
\put(8.9,7){\color{lgrey}\line(1,0){9.2}}
\put(34.9,7){\color{lgrey}\line(1,0){9.2}}
\put(7,20){\circle{3}}
\put(46,7){\circle{3}}
\put(7,18.1){\color{hgrey}\line(0,-1){9.2}}
\end{picture}
\caption{A picture illustrating the proof of Theorem~\ref{thm.poset-characterization}.\label{fig.2}}
\end{center}
\end{figure}

{\bf II.}\, The first step to verify the sufficiency is to observe that if $S$ is geodetically convex then,
for every $(u,v)\in\Cov (S)$ and $\permC\in S$, one has $u\prec_{\permC} v$.\\[0.4ex]
Following the definition of $\Cov (S)$ in Section~\ref{ssec.characterization} we
suppose specifically that an edge between enumerations $\permA\in S$ and $\permB\in \enuN\setminus S$ exists
with $\permA \stackrel{uv}{\longleftrightarrow} \permB$ and $u\prec_{\permA}v$, which implies $v\prec_{\permB}u$;
see Figure~\ref{fig.2}(a) for illustration.
Consider a geodesic $\mathfrak{g}$ between $\permA$ and $\permC$. Its nodes belong to $S$, by the
assumption that $S$ is geodetically convex.

Assume for a contradiction that $v\prec_{\permC} u$ giving $\{u,v\}\in\Inv [\permA,\permC]$ and,
by Lemma~\ref{lem.permut-graph}(ii), a unique edge of $\mathfrak{g}$ is labeled by $\{u,v\}$.
Let it be the one between $\permD\in S$ and $\permE\in S$, with $\permD$ closer to $\permA$.
Therefore, the section $\mathfrak{g}^{\prime}$ of $\mathfrak{g}$ between $\permA$ and $\permD$ is in $S_{u\prec v}$.
Because the mapping $\permG\mapsto\trans_{uv}\cdot\permG$ is an automorphism of the permutohedral graph
(see Section~\ref{ssec.permut-graph}), the section $\mathfrak{g}^{\prime}$ has a copy $\mathfrak{g}^{\prime\prime}$ in $S_{v\prec u}$, which is between $\permB$ and $\permE$\/; see Figure~\ref{fig.2}(b) for illustration.
The walk composed of the edge $[\permA,\permB]$ and $\mathfrak{g}^{\prime\prime}$
has the same length as the section of $\mathfrak{g}$ between $\permA$ and $\permE$ then.
Therefore, it is a geodesic between $\permA$ and $\permE$. By replacing the section of
$\mathfrak{g}$ between $\permA$ and $\permE$ by this one gets a geodesic between $\permA$
and $\permC$ containing $\permB$.
The assumption that $S$
is geodetically convex implies a contradictory conclusion that $\permB\in S$.
Thus, one necessarily has $u\prec_{\permC} v$.

{\bf III.}\, The next step is to show that, if $\emptyset\neq S\subseteq\enuN$ is a non-empty geodetically convex set
and $\permE\in\enuN\setminus S$ then there is $(u,v)\in\Cov (S)$ with
$v\prec_{\permE} u$.\\[0.4ex]
Consider a walk $\mathfrak{g}$ from $\permE$ to some $\permA\in S$ of
the least possible length among walks from $\permE$ to~$S$. By definition, $\mathfrak{g}$ is a geodesic between $\permE$
and $\permA$. Let $\permB$ be the last node of $\mathfrak{g}$ before $\permA$. Then necessarily
$\permB\in\enuN\setminus S$; let $\{u,v\}$ be the label of the edge between $\permB$ and $\permA$.
Without loss of generality assume $u\prec_{\permA} v$ (otherwise we exchange $u$ for $v$), which gives $(u,v)\in\Cov(S)$. By Lemma~\ref{lem.permut-graph}(ii) applied to $\mathfrak{g}$, $\{u,v\}$ is an inversion
between $\permE$ and $\permA$, which gives $v\prec_{\permE} u$.

{\bf IV.}\, As $|N|\geq 2$ is assumed, one has $\emptyset=(N\times N)^{\lGprec}$ and $\emptyset\in {\calX}^{\preced}$
(see Section~\ref{sec.linear-extensions}). Thus, to verify
the sufficiency it is enough to show $S\in {\calX}^{\preced}$
for a non-empty geodetically convex set $\emptyset\neq S\subseteq\enuN$.
Recall additionally from Section~\ref{sec.linear-extensions}
that, given $\permC\in\enuN$, $\permC\in\Cov (S)^{\lGprec}$ means
that $u\preceq_{\permC} v$ for any $(u,v)\in \Cov (S)$. Thus,
by Step II., one has $S\subseteq\Cov (S)^{\lGprec}$.
Analogously, by Step~III., one has $\enuN\setminus S ~\subseteq~ \enuN\setminus (\Cov(S)^{\lGprec})$, which gives
$\Cov (S)^{\lGprec}\subseteq S$. We have thus shown $S=\Cov (S)^{\lGprec}$, which means
$S\in {\calX}^{\preced}$ (see Section~\ref{ssec.app.Galois}).

{\bf V.}\, The second claim in Theorem~\ref{thm.poset-characterization} then follows from its first claim
and Lemma~\ref{lem.poset-lattice}(ii).
\end{proof}

\section{Proof of Theorem~\ref{thm.poset-to-extensions}}\label{app.thm.poset-to-extensions}

The notation introduced in Sections~\ref{ssec.trans-closure}, \ref{ssec.face-lattice}, \ref{sec.linear-extensions}, and \ref{ssec.characterization} is utilized in the proof.

\noindent{\bfseries Theorem~\ref{thm.poset-to-extensions}:} ~~\rm
Let $P$ be a poset over $N$ and $S={\cal L}(P)$. Then, for any pair $u,v\in N$, one has
\begin{itemize}
\item[(i)] $u\|v\,[P]$ if and only if $\{u,v\}\in \Inv (S)$,
\item[(ii)] $u\cover v\,\,[P]$ if and only if  $(u,v)\in\Cov (S)$,
\item[(iii)] $(u,v)\in P\setminus\diag$ if and only if  $(u,v)\in\tr(\Cov (S))$,
\item[(iv)]  $(u,v)\in P$ if and only if  $(u,v)\in\tr(\diag\cup\Cov (S))$.
\end{itemize}

\begin{proof}
We first verify a series of preparatory observations for a given pair $u,v\in N$ of variables.\\[0.4ex]
{\bf I.}\, If either $u\|v\,[P]$ or $u\cover v\,\,[P]$ then $\exists\, \permA\in {\cal L}(P) \,:\ u\covA v$\,.\\[0.4ex]
We put $A:=\{\, a\in N\setminus v\,:\ (a,u)\in P ~\mbox{or}~ (a,v)\in P\,\}$, $B:=N\setminus A$
and derive the next 3 facts.
\begin{itemize}
\item The variable $v$ is initial in $B$, that is, $v\in B$ and
$\neg\,[\,\exists\, b\in B\setminus v ~:\ (b,v)\in P\,]$.
\end{itemize}
Indeed, otherwise, by the definition of $A$, one has $b\in A$, contradicting the definition of $B$.
\begin{itemize}
\item The variable $u$ is terminal in $A$, that is, $u\in A$ and
$\neg\,[\,\exists\, a\in A\setminus u ~:\ (u,a)\in P\,]$.
\end{itemize}
Assume for a contradiction that such $a\neq u$ exists. Then $a\in A$ means $a\neq v$ and either $(a,u)\in P$
or $(a,v)\in P$. As $P$ is a poset, the first option $(a,u),(u,a)\in P$ contradicts $a\neq u$. The second option $(u,a),(a,v)\in P$ contradicts both $u\|v\,[P]$ (by transitivity of $P$) and $u\cover v\,\,[P]$.
\begin{itemize}
\item The set $A$ precedes the set $B$ in $P$, that is,
$\neg\,[\,\exists\, b\in B ~\&~ a\in A ~:\ (b,a)\in P\,]$.
\end{itemize}
Assume for a contradiction that such $b$ and $a$ exist and distinguish the cases $b\neq v$ and $b=v$.
If $b\neq v$ then, by transitivity of $P$, $b\in A$, contradicting the definition of $B$. If $b=v$ then
the option $(v,a),(a,u)\in P$ implies, by transitivity of $P$, $(v,u)\in P$ and contradicts both
$u\|v\,[P]$ and $u\cover v\,\,[P]$. The other option $(v,a),(a,v)\in P$ contradicts $a\neq v$, as $P$ is a poset.\\[0.4ex]
Consider an ordered partition $|A_{1}|A_{2}|A_{3}|A_{4}| \,:=\, |A\setminus u\,|\,u\,|\,v\,|\,B\setminus v|$ of $N$.
The above derived facts altogether say that $P$ is consistent with it:
$(c,d)\in P,\, c\in A_{i},\, d\in A_{j}  ~\Rightarrow~  i\leq j$. As any
relation $P_{i}:=P\cap (A_{i}\times A_{i})$ is a poset on $A_{i}$, for every $i=1,2,3,4$,
there exists $\permA(i)\in {\cal L}(P_{i})$. One can compose an enumeration $\permA\in\enuN$
of these: $\permA \,:=\, |\permA(1)|\permA(2)|\permA(3)|\permA(4)|$. It follows from the construction
that $\permA\in {\cal L}(P)$ and also $u\covA v$.\\[0.6ex]
{\bf II.}\, If $u\|v\,[P] ~\&~ \permA\in {\cal L}(P) \,:\ u\covA v$
then $\permB\in {\cal L}(P)$ for $\permB\in\enuN$ specified by $\permA \stackrel{uv}{\longleftrightarrow} \permB$.\\[0.4ex]
Having $(c,d)\in P$, the assumption $u\|v\,[P]$ implies $(u,v)\neq (c,d)\neq (v,u)$ and $\permA\in {\cal L}(P)$
implies $c\preceq_{\permA}d$. Now, $\permA \stackrel{uv}{\longleftrightarrow} \permB$ says that $\permA$ and
$\permB$ only differ in the mutual position of $u$ and $v$, implying thus $c\preceq_{\permB}d$.
Hence, we have verified $\permB\in {\cal L}(P)$.
\\[0.6ex]
{\bf III.}\, If $\exists\, \permA,\permB\in {\cal L}(P) ~:\ \permA \stackrel{uv}{\longleftrightarrow} \permB$
then $u\|v\,[P]$.\\[0.4ex]
Assume specifically $u\prec_{\permA} v$ and  $v\prec_{\permB} u$, which necessitates
$\neg[ v\preceq_{\permA} u]$ and $\neg[ u\preceq_{\permB} v]$. The option $(v,u)\in P$ thus
forces $\permA\not\in {\cal L}(P)$, yielding a contradiction. The option $(u,v)\in P$
enforces $\permB\not\in {\cal L}(P)$, giving a contradiction, too. So, we have verified $u\|v\,[P]$.\\[0.6ex]
{\bf IV.}\, If $u\cover v\,\,[P] ~\&~ \permA\in {\cal L}(P) \,:\ u\covA v$
then $\permB\not\in {\cal L}(P)$ for $\permB\in\enuN$ specified by $\permA \stackrel{uv}{\longleftrightarrow} \permB$.\\[0.4ex]
Assume $\permB\in {\cal L}(P)$ for a contradiction. Then there are $\permA,\permB\in {\cal L}(P)$ with
$\permA \stackrel{uv}{\longleftrightarrow} \permB$ and, by Step\,III., $u\|v\,[P]$, yielding a contradiction
with the assumption $u\cover v\,\,[P]$.\\[0.6ex]
{\bf V.}\, If $\exists\, \permA\in {\cal L}(P) ~\&~ \permB\in \enuN\setminus {\cal L}(P) ~:\ \permA \stackrel{uv}{\longleftrightarrow} \permB ~\&~ u\prec_{\permA} v$ then $u\cover v\,\,[P]$.\\[0.4ex]
Realize that $\permA \stackrel{uv}{\longleftrightarrow} \permB$ and $u\prec_{\permA} v$ imply $u\covA v$.
If $u\|v\,[P]$ then Step~II.\/ implies a contradictory
conclusion $\permB\in {\cal L}(P)$. If $(v,u)\in P$ then $\permA\in {\cal L}(P)$ gives  $v\preceq_{\permA}u$
contradicting $u\covA v$. Thus, necessarily $(u,v)\in P$. Assume for a contradiction that $w\in N\setminus uv$ exists
with $(u,w),(w,v)\in P$. Then $\permA\in {\cal L}(P)$ implies $u\prec_{\permA} w\prec_{\permA} v$ contradicting
$u\covA v$. Therefore, $u\cover v\,\,[P]$.
\medskip

Now, the claim (i) is that $u\|v\,[P]$ if and only if $\exists\, \permA,\permB\in {\cal L}(P) \,:\ \permA \stackrel{uv}{\longleftrightarrow} \permB$, by the definition of $\Inv(S)$. The necessity here follows from Steps~I.\/ and II., the sufficiency from Step~III.
\smallskip

The claim (ii) means, by the definition of $\Cov(S)$, that $u\cover v\,\,[P]$ if and only if one has $\exists\, \permA\in {\cal L}(P) ~\&~ \permB\in \enuN\setminus {\cal L}(P) ~:\ \permA \stackrel{uv}{\longleftrightarrow} \permB ~\&~ u\prec_{\permA} v$. The necessity follows from Steps~I.\/ and IV., the sufficiency from Step~V. Further auxiliary observation is needed to derive the rest.\\[0.6ex]
{\bf VI.}\, The relation $\Cov ({\cal L}(P))$ on $N$ is acyclic.\\[0.4ex]
It is enough to verify $\Cov ({\cal L}(P))\subseteq P\setminus\diag$, because $P\setminus\diag$ is acyclic (see Fact~\ref{fct.posets} in Section~\ref{ssec.trans-closure}). Given $(u,v)\in \Cov ({\cal L}(P))$, the definition of $\Cov(S)$ in Section~\ref{ssec.characterization} implies the existence
of  $\permA\in {\cal L}(P)$ and  $\permB\in \enuN\setminus {\cal L}(P)$ such that $\permA \stackrel{uv}{\longleftrightarrow} \permB ~\&~ u\prec_{\permA} v$.
Step~V.\/ gives $u\cover v\,\,[P]$, and, thus, $(u,v)\in P\setminus\diag$.
\medskip

For the necessity in (iii) consider $(u,v)\in P\setminus\diag$. Because $N$ is finite, one can
show by repeated ``including" elements $w\in N$ between $u$ and $v$, that $\exists\, u=u_{1},\ldots,u_{k}=v$,
$k\geq 2$ such that $u_{i}\cover u_{i+1}\,\,[P]$ for $i=1,\ldots,k-1$. Then the application
of (ii) gives $(u_{i},u_{i+1})\in\Cov (S)$ for $i=1,\ldots,k-1$, implying
$(u,v)\in\tr(\Cov (S))$.
\smallskip

For the sufficiency in (iii) consider $(u,v)\in\tr(\Cov (S))$ and $u=u_{1},\ldots,u_{k}=v$,
$k\geq 2$ such that $(u_{i},u_{i+1})\in\Cov (S)$ for $i=1,\ldots,k-1$. By (ii), one has
$u_{i}\cover u_{i+1}\,\,[P]$ for these~$i$, and the transitivity of $P$ yields $(u,v)\in P$.
Assume for a contradiction $u=v$, in which case one has $\exists\, u=u_{1},\ldots,u_{k}=u$,
$k\geq 2$ with $(u_{i},u_{i+1})\in\Cov (S)$ for $i=1,\ldots,k-1$, meaning that
$\Cov({\cal L}(P))$ contains a cycle (see Section~\ref{ssec.trans-closure}), yielding
a contradiction, by Step~VI.
\smallskip

The claim (iv) follows directly from (iii): as $P=\diag\cup (P\setminus\diag)$, one
can apply evident equality $\diag\cup\tr(R)=\tr(\diag\cup R)$
valid for any $R\subseteq N\times N$.
\end{proof}

\section{Proof of Lemma~\ref{lem.encoding-posets}}\label{app.lem.encoding-posets}

\noindent{\bfseries Lemma~\ref{lem.encoding-posets}:} ~~\rm
Let $\permD\in\enuN$ be a fixed reference enumeration of $N$, where $|N|\geq 2$.
\begin{itemize}
\item[(i)] Given a poset $P$ on $N$, one has ${\cal L}(P)={\cal L}^{\permD}[A,B]$, where
$$
{\cal L}^{\permD}[A,B] \,:=\, \{\, \permC\in\enuN\,:\  A\subseteq T_{\permD}\setminus T_{\permC}\subseteq B\,\}
$$
with $A:=(T_{\permD}\cap P^{op})\setminus\diag$ and $B:=T_{\permD}\setminus P$, clearly
satisfying $A\subseteq B\subseteq T_{\permD}\setminus\diag$.
\item[(ii)] A set $S\subseteq\enuN$ of enumerations of $N$ has the form ${\cal L}^{\permD}[A,B]$ for
a pair $[A,B]$ of subsets of\/ $T_{\permD}\setminus\diag$ if and only if it is geodetically convex in the
permutohedral graph.
\item[(iii)] Given a non-empty geodetically convex set $S\subseteq\enuN$, the pair $[A,B]$
with $S={\cal L}^{\permD}[A,B]$ and minimal set difference $B\setminus A$ is unique.
It has the forms $[\,\,\bigcap_{\permC\in S} T_{\permD}\setminus T_{\permC}, \bigcup_{\permC\in S} T_{\permD}\setminus T_{\permC}\,]$ and $[\,(T_{\permD}\cap P^{op})\setminus\diag,T_{\permD}\setminus P\,]$, where $P$ is the unique poset on $N$ with $S={\cal L}(P)$.
\end{itemize}

\begin{proof}
Throughout the proof we denote\, $T:=T_{\permD}\setminus\diag$.

To prove (i) realize that ${\cal L}(P)= \{\, \permC\in\enuN\,:\  P\subseteq T_{\permC}\,\}$. Clearly, $P\subseteq T_{\permC}$ is equivalent to the inclusion
$(N\times N)\setminus T_{\permC}\subseteq (N\times N)\setminus P$.
Because $N\times N$ disjunctively decomposes as $T\cup\diag\cup T^{op}$,
this inclusion breaks into two conditions, namely into $T\setminus T_{\permC}\subseteq T\setminus P=B$ and\, $T^{op}\setminus T_{\permC}\subseteq T^{op}\setminus P$.
Because of one-to-one correspondence $R\longleftrightarrow R^{op}$, the latter condition has an equivalent writing\, $T\setminus T_{\permC}^{op}\subseteq T\setminus P^{op}$.
However, as $T\setminus T_{\permC}^{op}=T\cap T_{\permC}$, it means $T\cap T_{\permC}\cap P^{op}=\emptyset$, that is, $A=T\cap P^{op}\subseteq T\setminus  T_{\permC}$. As $P$ is a poset, $P\cap P^{op}=\diag$ yields $A\subseteq B$.

For the necessity in (ii) assume $\permA,\permB\in {\cal L}^{\permD}[A,B]$ and $\permC\in\enuN$ between them.
In case $(a,b)\in A\subseteq T$ one has $(b,a)\in T_{\permA}\cap T_{\permB}$ and,
by Lemma~\ref{lem.permut-graph}(iii), $(b,a)\in T_{\permC}$, implying $(a,b)\in T\setminus T_{\permC}$.
The inclusion $T\setminus T_{\permC}\subseteq B$ can be shown in the form $T\setminus B\subseteq T_{\permC}$:
if $(a,b)\in T\setminus B$ then $(a,b)\in T_{\permA}\cap T_{\permB}$, and, by Lemma~\ref{lem.permut-graph}(iii), $(a,b)\in T_{\permC}$. Therefore, $\permC\in {\cal L}^{\permD}[A,B]$, which was the goal.

For the sufficiency in (ii) assume $S\neq\emptyset$, as otherwise $S=\emptyset={\cal L}^{\permD}[T,\emptyset]$.
Apply Theorem~\ref{thm.poset-characterization} to $S$ to get
the unique poset $P$ on $N$ with $S={\cal L}(P)$ and then use (i).

To prove (iii), let us put $A(S):=\bigcap_{\permC\in S} T\setminus T_{\permC}$ and $B(S):=\bigcup_{\permC\in S} T\setminus T_{\permC}$. We known already from the above arguments concerning (ii) that $S={\cal L}(P)$ for a poset on $N$ and $S={\cal L}^{\permD}[A,B]$ for some $[A,B]$. The inclusions $A\subseteq A(S)\subseteq T\setminus T_{\permC}\subseteq B(S)\subseteq B$ valid for every $\permC\in S$ allow one to show $S={\cal L}^{\permD}[A(S),B(S)]$, which yields the first formula. To get the second formula one needs to show both $A(S)\subseteq T\cap P^{op}$ and $(T\setminus P)\subseteq B(S)$ and combine these inclusions with (i).

In case of $(u,v)\in A(S)=\bigcap_{\permC\in S} T\setminus T_{\permC}$ one cannot have $(u,v)\in P$
because it leads to a contradictory conclusion $(u,v)\in T_{\permC}$ for arbitrary $\permC\in {\cal L}(P)=S$.
Additionally, $u\|v\,[P]$ implies, by Theorem~\ref{thm.poset-to-extensions}(i), that there are
$\permA,\permB\in S={\cal L}(P)$ with $\permA \stackrel{uv}{\longleftrightarrow} \permB$, implying
that either $(u,v)\in T_{\permA}$ or  $(u,v)\in T_{\permB}$, both contradicting $(u,v)\in A(S)$.
Hence, necessarily $(u,v)\in P^{op}$.

In case $(u,v)\in T\setminus P$ one either has $(v,u)\in P$ or $u\|v\,[P]$. In the former case one
has $(v,u)\in T_{\permC}$ for arbitrary $\permC\in {\cal L}(P)=S$, yielding $(u,v)\not\in T_{\permC}$ and $(u,v)\in B(S)$.
In the latter case, Theorem~\ref{thm.poset-to-extensions}(i) implies that there are
$\permA,\permB\in S={\cal L}(P)$ with $\permA \stackrel{uv}{\longleftrightarrow} \permB$,
which gives either $(u,v)\not\in T_{\permA}$ or $(u,v)\not\in T_{\permB}$, with both options
yielding $(u,v)\in B(S)$.
\end{proof}

\section{Proof of Lemma~\ref{lem.braid-lattice}}\label{app.lem.braid-lattice}

\noindent{\bfseries Lemma~\ref{lem.braid-lattice}:} ~~\rm
Let $({\calX}^{\braid},\subseteq )$ and $({\calY}^{\braid},\subseteq )$ denote the lattices defined using the
Galois connections $\lGbraid$ and $\rGbraid$ based on the incidence relation $\braid$ introduced in Section~\ref{sssec.preposet-cones}.
\begin{itemize}
\item[(i)] Then $T\in {\calY}^{\braid}$ if and only if\/ $T$ is a preposet on $N$.
\item[(ii)] A preposet $T\subseteq N\times N$ is a poset if and only if its braid cone $\nor_{T}$ has
the full dimension $n$.
\item[(iii)] Given a poset $P$ on $N$, one has ${\cal L}(P)=\{\, \permA\in\enuN \,:\ \nor^{\permA}\subseteq\nor_{P}\,\}$
and $\nor_{P}=\bigcup_{\permA\in {\cal L}(P)} \nor^{\permA}$.
\end{itemize}

\begin{proof}
For the necessity in (i) assume $T=\nor^{\rGbraid}=\{\,(u,v)\in N\times N\, :\ \forall\, x\in\nor\quad x_{u}\leq x_{v} \,\}$ for some $\nor\subseteq {\dv R}^{N}$ (see Section~\ref{ssec.app.Galois}). Evidently, $T$ is both reflexive and transitive.

For the sufficiency in (i) assume that $T$ is a preposet on $N$ and denote by ${\cal E}$ the set of equivalence classes of $E:=T\cap T^{op}$. Put $\nor:=\nor_{T}=T^{\lGbraid}$; it is immediate
that\/ $T\subseteq \nor^{\rGbraid}$. To verify the other inclusion $\nor^{\rGbraid}\subseteq T$, we will show
$(a,b)\in (N\times N)\setminus T ~\Rightarrow~ (a,b)\not\in \nor^{\rGbraid}$.\\[0.4ex]
For this purpose, we interpret $T$ in the form of a partial order $\preceq$ on ${\cal E}$ (see Section~\ref{ssec.posets-defintions}).
This poset $({\cal E},\preceq)$ can additionally be viewed as a transitive directed acyclic graph $G$ over ${\cal E}$ (Section~\ref{ssec.poset-as-DAG}).
The assumption $(a,b)\not\in T$ thus implies that there are $A,B\in {\cal E}$ with $a\in A$, $b\in B$ and
$\neg[A\preceq B$], that is, both $A\neq B$ and $\neg[\,A\to B ~\mbox{in $G$}\,]$. Therefore, an ordering $U_{1},\ldots, U_{m}$, $m\geq 1$, of elements of ${\cal E}$ exists, which is in concordance with $G$
and in which $B$ strictly precedes $A$.
Define $x\in {\dv R}^{N}$ as the ``rank vector" of this block-enumeration:
put, for any $u\in N$, $x_{u}:=i$ where $u\in U_{i}$. Then $(u,v)\in T$ implies either $[\exists\,i\,: u,v\in U_{i}]$ or $[\exists\,i<j\,: u\in U_{i}\,\&\,v\in U_{j}]$, yielding $x_{u}\leq x_{v}$.
Hence, $x\in \nor_{T}=\nor$ but $x_{b}<x_{a}$, and $(a,b)\not\in \nor^{\rGbraid}$, which was the goal.
\smallskip

For the necessity in (ii) realize that ${\cal L}(T)\neq\emptyset$. Thus, $T\subseteq T_{\permA}$ for some
$\permA\in {\cal L}(T)\subseteq\enuN$ and the anti-isomorphism of $({\calX}^{\braid},\subseteq )$ and $({\calY}^{\braid},\subseteq )$ yields $\nor^{\permA}\subseteq  \nor_{T}$.

For the sufficiency in (ii) observe that there exists $y\in\nor_{T}$ whose all components are distinct.
Because of $(u,v)\in T \,\Rightarrow\, y_{u}\leq y_{v}$ one has
$(u,v),(v,u)\in T \,\Rightarrow\, y_{u}=y_{v}$, enforcing $u=v$.

To prove (iii) realize that $\permA\in {\cal L}(P) \,\Leftrightarrow\, P\subseteq T_{\permA} \,\Leftrightarrow\,
\nor^{\permA}\subseteq\nor_{P}$, owing to the properties of Galois connections. Hence,
$\bigcup_{\permA\in {\cal L}(P)} \nor^{\permA}\subseteq \nor_{P}$. For the other inclusion
realize that, if $y\in \nor_{P}$ has all components distinct then unique $\permB\in\enuN$ exists with $y\in \nor^{\permB}$.
The fact $y_{u}<y_{v} \,\Leftrightarrow\, u\prec_{\permB} v$ then implies $\permB\in {\cal L}(P)$.
The fact that the vectors $y\in {\dv R}^{N}$ with distinct components are dense in $\nor_{P}$ allows one to
derive the other inclusion by applying the closure in Euclidean topology.
\end{proof}

\section{Proof of Lemma~\ref{lem.topol-lattice}}\label{app.lem.topol-lattice}

\noindent{\bfseries Lemma~\ref{lem.topol-lattice}:} ~~\rm
Let $({\calX}^{\topol},\subseteq )$ and $({\calY}^{\topol},\subseteq )$ denote the lattices defined using the
Galois connections $\lGtopol$ and $\rGtopol$ based on the incidence relation $\topol$ defined in Section~\ref{sssec.preposet-topologies}.
\begin{itemize}
\item[(i)] Then $T\in {\calY}^{\topol}$ if and only if\/ $T$ is a preposet on $N$.
\item[(ii)] A set system $\calD\subseteq\caP$ belongs to the lattice ${\calX}^{\topol}$ if and only if
both $\emptyset,N\in\calD$ and the implication $D,E\in\calD \,\Rightarrow\, D\cap E,D\cup E\in\calD$
holds (= $\calD$ is a topology on $N$).
\item[(iii)] A preposet $T\subseteq N\times N$ is a poset if and only if
the topology $\calD_{T}$ distinguishes points: for any distinct $u,v\in N$
there is $D\in\calD_{T}$ such that either $[\,u\in D ~\&~ v\not\in D \,]$ or
$[\,v\in D ~\&~ u\not\in D \,]$. This happens exactly when there is $\permA\in\enuN$ such that
$\chain_{\permA}\subseteq\calD_{T}$.
\item[(iv)] Given a braid cone $\nor\subseteq {\dv R}^{N}$, the set system
$\calD_{\nor} :=\{\, D\subseteq N\, :\ -\chi_{D}\in\nor\,\}$ is the corresponding topology.
Given a topology $\calD$, the cone
$\nor_{\calD} := \cone(\{\chi_{N}\}\cup\{-\chi_{D}:D\in\calD\})$,
where $\cone(\nor)$ denotes the conic hull of $\nor\subseteq {\dv R}^{N}$, is the corresponding braid cone.
\item[(v)] Given a poset $P$ on $N$, one has ${\cal L}(P)=\{\, \permA\in\enuN \,:\ \chain_{\permA}\subseteq\calD_{P}\,\}$
and $\calD_{P}=\bigcup_{\permA\in {\cal L}(P)} \chain_{\permA}$.
\end{itemize}

\begin{proof}
For the necessity in (i) assume\, $T=\calD^{\rGtopol}:=\{\,(u,v)\in N\times N\, :\ \forall\, D\in\calD~~\,
v\in D \Rightarrow u\in D \,\}$ for some $\calD\subseteq\caP$ (see Section~\ref{ssec.app.Galois}). Evidently, $T$ is both reflexive and transitive.

For the sufficiency in (i) assume that $T$ is a preposet on $N$; it is immediate that $T\subseteq (\calD_{T})^{\rGtopol}$. For any $v\in N$ we put $I(v):=\{\, w\in N\,:\ (w,v)\in T\,\}$; the reflexivity of\/ $T$ implies $v\in I(v)$, while the transitivity of\/ $T$ allows one to observe $I(v)\in \calD_{T}$.
Given $(u,v)\in (\calD_{T})^{\rGtopol}$, these two facts about $I(v)$ and the definition of $\calD^{\rGtopol}$ imply $u\in I(v)$, which means $(u,v)\in T$. Therefore, $(\calD_{T})^{\rGtopol}\subseteq T$ and $T=(\calD_{T})^{\rGtopol}$ yields $T\in {\calY}^{\topol}$.

The necessity in (ii) is immediate: assume $T\subseteq N\times N$ and realize that the definition of\/ $\calD_{T}$ allows one to verify easily the conditions in (ii) for $\calD_{T}$.

For the sufficiency in (ii) assume that $\calD\subseteq\caP$ is a topology on $N$ and put $T:=\calD^{\rGtopol}$;
it is immediate that $\calD\subseteq T^{\lGtopol}=\calD_{T}$. For any $v\in N$ we put $D(v):=\bigcap_{\,D\in\calD\,:\,v\in D} D$.
Then $N\in\calD$ implies $v\in D(v)$ and, since $\calD$ is closed under intersection, $D(v)\in\calD$.
Thus, $E\subseteq\bigcup_{v\in E} D(v)$ for any $E\subseteq N$. Given $E\in\calD_{T}$, the goal is to verify
the other inclusion $\bigcup_{v\in E} D(v)\subseteq E$.\\[0.4ex]
Indeed, having $w\in \bigcup_{v\in E} D(v)$, there exists $v\in E$ with $w\in D(v)$. The latter condition means
$[\,\forall\, D\in\calD~~ v\in D \Rightarrow w\in D\,]$, which is nothing but $(w,v)\in \calD^{\rGtopol}=T$.
By the definition of $\calD_{T}$, the facts $E\in\calD_{T}$, $(w,v)\in T$, and $v\in E$ imply
$w\in E$, concluding the proof of the inclusion.\\[0.4ex]
Therefore, since $\emptyset\in\calD$ and $\calD$ is closed under union, one gets $E=\bigcup_{v\in E} D(v)\in\calD$.
This concludes the proof of $E\in \calD_{T}\,\Rightarrow\, E\in\calD$ and $\calD=\calD_{T}=T^{\lGtopol}$ yields $\calD\in {\calX}^{\topol}$.

For (iii) note that if $T$ is poset then $T\subseteq T_{\permA}$ for some enumeration $\permA\in\enuN$ and the anti-isomorphism of $({\calX}^{\topol},\subseteq )$ and $({\calY}^{\topol},\subseteq )$ yields $\chain_{\permA}\subseteq\calD_{T}$. Hence, $\calD_{T}$ distinguishes points.
Conversely, if $\calD_{T}$ distinguishes points then there is no pair of distinct $u,v\in N$ with $(u,v),(v,u)\in T$, for otherwise $\forall\, D\in\calD_{T}~\, u\in D \Leftrightarrow v\in D$.
Thus, $T\cap T^{op}=\diag$ and\/ $T$ is anti-symmetric.

To prove (iv), assume $\nor=\nor_{T}$ and $\calD=\calD_{T}$ for a preposet $T$ on~$N$. Given $D\subseteq N$ and $(u,v)\in N\times N$, re-write $D\,\topol\, (u,v)$, that is, $[\,v\in D \Rightarrow u\in D\,]$ as $[\, (\chi_{D})_{v}=1 \,\Rightarrow\, (\chi_{D})_{u}=1\,]$ and then as $(-\chi_{D})_{u}\leq (-\chi_{D})_{v}$, which means $(-\chi_{D})\braid (u,v)$. This allows one to observe that $D\in\calD_{T}$ if and only if $-\chi_{D}\in\nor_{T}$, implying the first claim in (iv).

The above observation also gives the inclusion $\cone(\{\chi_{N}\}\cup\{-\chi_{D}:D\in\calD_{T}\})\subseteq \nor_{T}$, because $\nor_{T}$ is a cone. Therefore, the verification of the other inclusion
$$
\nor_{T}\subseteq \cone(\{\chi_{N}\}\cup\{-\chi_{D}:D\in\calD_{T}\})
$$
is needed to complete the proof of the second claim in (iv).\\[0.3ex]
For this purpose, consider the set ${\cal E}$ of equivalence classes of $E:=T\cap T^{op}$ and view $T$ as a partial order on ${\cal E}$ and, thus, as a transitive directed acyclic graph $G$ over ${\cal E}$ (see Sections~\ref{ssec.posets-defintions} and \ref{ssec.poset-as-DAG}).
Observe that any vector $x\in\nor_{T}$ complies with $G$ in the following sense:
if $u\in U\in {\cal E}$, $v\in V\in {\cal E}$ and $U=V$ or $U\to V$ in $G$ then $x_{u}\leq x_{v}$. Hence, $x$ is ``constant" on sets from
${\cal E}$, that is, $x_{u}=x_{v}$ if $(u,v)\in E$.
This allows one to construct an ordering $U(1),\ldots, U(m)$, $m\geq 1$, of (all) elements of ${\cal E}$, which
is both in concordance with $G$ and respects $x$ in the sense that
$$
u\in U(i),~ v\in U(j),~ \mbox{and}~ i\leq j\quad \mbox{implies}\quad x_{u}\leq x_{v}
\,.
$$
Let us introduce $x(i)$ for $i=1,\ldots, m$ as the shared value $x_{u}$ for $u\in U(i)$ and write
$$
x= x(m)\cdot\chi_{N} \,+\, \sum_{i=1}^{m-1}\, (x(i+1)- x(i))\cdot (-\chi_{\,U(1)\ldots U(i)})\,,
$$
which equality can be verified by substituting any $u\in U(j)$ for $j\in [m]$.
Since $x(i+1)-x(i)\geq 0$ for $i\in[m-1]$, the vector $x$ belongs to the conic hull of $\{\,-\chi_{\,U(1)\ldots U(i)}\,:\ i\in [m-1]\,\}\cup\{\pm\chi_{N}\}$.
The concordance of $U(1),\ldots, U(m)$ with $G$ implies that $U(1)\ldots U(i)\in\calD_{T}$ for $i\in [m]$, which thus yields the desired inclusion and completes the proof of the second claim in (iv).

To prove (v) realize that $\permA\in {\cal L}(P) \,\Leftrightarrow\, P\subseteq T_{\permA} \,\Leftrightarrow\,
\chain_{\permA}\subseteq\calD_{P}$, owing to the properties of Galois connections. Hence,
$\bigcup_{\permA\in {\cal L}(P)} \chain_{\permA}\subseteq \calD_{P}$. For the other inclusion
take $A\in\calD_{P}$, put $B:=N\setminus A$, and observe there is no $(b,a)\in P$ with $b\in B$ and $a\in A$.
In case $A=\emptyset$ or $A=N$ one has $A\in\chain_{\permA}$ for any $\permA\in {\cal L}(P)$.
As $P_{A}:=P\cap (A\times A)$ is a poset on $A$, there is an enumeration $\permA(A)\in {\cal L}(P_{A})$, analogously
with $B$, and concatenation $\permA:=|\permA(A)|\permA(B)|$ yields $\permA\in {\cal L}(P)$ with $A\in\chain_{\permA}$.
\end{proof}

\section{Essay: reasoning for our notational conventions}\label{app.remark-conventions}
This long remark (on a request of a reviewer) is to explain in detail our reasons for adopting the conventions from Section~\ref{ssec.notation-convention}, namely
that the ground set $N$ is regarded as an unordered set, while the discrete interval $[n]:=\{1,\ldots,n\}$, where $n=|N|$, is regarded as a totally
ordered set, with the ordering being the natural one, inherited from the set ${\dv N}$ of natural numbers. An immediate consequence of these assumptions
is the fact $N\neq [n]$, which is a crucial point.

The main motive for our approach was to come with a consistent and clear presentation of our results and proofs. This is because in some former papers addressing related topics \cite{PRW08,MPSSW09,Cri11,MUWY18} a complete opposite convention was accepted from the very beginning: the ground set was simply identified with the discrete interval by putting $N:=[n]$. When studying these works we encountered discrepancies and ambiguous notation which, we believe,
stem precisely from that differing convention of $N:=[n]$. To avoid similar inconsistencies in our own paper we came with our conventions.
Let us point at an elementary difference between the conventions.
\smallskip

In our paper we deal with bijections between $N$ and $[n]$ and these bijections
encode (different) total orders on the ground set $N$. Specifically, given a {\em rank vector\/} $\rankv:N\to [n]$,
the numerical value $i\in [n]$ assigned to an element $u\in N$ of the ground set by $\rankv$ encodes the position of $u$ in the total order and, conversely, the inverse bijection $\rankv_{-1}$ assigns an element $u$ of the ground set to the $i$-th position.
In our approach, when $N \neq [n]$, the inverse of such a rank vector $\rankv:N\to [n]$ can never be confused with another such a rank vector $\varrankv:N\to [n]$, because the inverse mapping to $\rankv$ has a different domain $[n]$: note that $[n]\neq N$ by our convention.
The same principle holds if we decide to encode total orders on $N$ alternatively, namely by the
inverse bijections from $[n]$ onto $N$, which we name {\em enumerations}.
\smallskip

The papers adopting the competing convention also utilized bijections to encode total orders on the ground set $N:=[n]$. However, these are, in fact, self-bijections of the ground set~$[n]$.
Therefore, suddenly, a strange ambiguity occurs when one is trying to interpret a self-bijection $\permD:[n]\to [n]$ as a total order on $[n]$. Indeed, both sides of the transformation $\permD$, its domain and its image, can be interpreted both as the ground set and as the set of numerical values. These two sets, the domain and the image, are indistinguishable under this convention and $\permD$ can be interpreted as a total order on $[n]$ in {\bf two different ways}, depending on whether the domain or the image is interpreted as the ground set. For example, if $N:=[4]=\{1,2,3,4\}$ and $\permD$ is given by $\permD(1)=2$, $\permD(2)=4$, $\permD(3)=3$, and $\permD(4)=1$, then the former interpretation leads to the total order $4<1<3<2$ while the latter interpretation leads to $2<4<3<1$. That means, an additional convention is needed to interpret a self-bijection $\permD$ as a total order: either the domain or the image of $\permD$ is interpreted as the ground set. This is just an elementary nuisance but even worse interpretation/ambiguity problems occur with advanced tasks.
\medskip

One of the reviewers requested a concrete example to clarify why we regard our ground set as unordered and why we distinguish it from the discrete interval $[n]$. We believe that perhaps the best example illustrating both the utility of our convention and the difference from the competing approach is the problem of characterizing the group of structure-preserving self-transformations (= the symmetry group) of the set of all total orders on a given ground set. This problem is equivalent to describing the symmetry group of the permutohedron, a task recently resolved in \cite{Stu26}, where a previously unproven statement from \cite{Cri11} was confirmed.
\medskip

Let us first describe to the reader how, under {\bf our convention}, things become transparent and elements mutually distinguishable.
Our approach allows for a clear distinction between the {\em combinatorial\/} and {\em geometric\/} ways of encoding total orders on the ground set $N$. The combinatorial way uses the bijections $\permA:[n]\to N$, called {\em enumerations\/} by us, while the geometric way uses the bijections $\rankv:N\to [n]$, called {\em rank vectors\/} by us. Observe that the set of total orders on the ground set {\em does not form a group\/} in our approach. This is because all total orders on this unordered set are equally valid; consequently, there is no distinguished element within this collection of total orders like the identity element within a group.

On the other hand, the geometric view allows one to interpret the set of total orders on our unordered ground set $N$
as the set of vertices of a particular polytope $\perN$ in ${\dv R}^{N}$, called the {\em permutohedron} by us,
and this geometric object {\em admits group actions\/} from the symmetric group on an $n$-element set, $n=|N|$.
The point is that \underline{two} \underline{different} group actions on the permutohedron are possible, and our
notational convention allows one to discern them clearly.

Of course, we are interested in group actions on\/ $\perN$ that result in structure-preserving self-transformations of\/ $\perN$,
that is, we wish to describe the {\em symmetry group\/} of\/ $\perN$. The prominent group actions are those from the group $\sym_{N}$
of permutations of our ground set $N$ because these always constitute automorphisms of\/ $\perN$. These morphisms are
associated with the elements of $\sym_{N}$, that is, with self-bijections $\alpha:N\to N$ (of $N$).
Observe that there is no need for an additional convention regarding whether the action of an element of\/ $\alpha\in\sym_{N}$
to a rank vector $\rankv$ is defined as the composition with $\rankv$ either from the left of $\rankv$ or from the right of~$\rankv$: since $N\neq [n]$, there is only one option.
The same applies to the alternative mode in which the total orders on $N$ are encoded by enumerations of $N$.

One can also consider actions from the group $\sym_{[n]}$ of permutations of the discrete interval~$[n]$, which are associated with self-bijections $\omega:[n]\to [n]$. However, if $n\geq 3$ then not all actions from $\sym_{[n]}$ constitute automorphisms of\/ $\perN$.
The only exception, besides the identity on $[n]$, is the ``reversing" mapping $\omega(i):= n+1-i$ for $i\in [n]$, whose respective
automorphism is the central symmetry of\/ $\perN$. By \cite{Stu26}, if $n\geq 3$, the symmetry group of\/ $\perN$ is the direct product of $\sym_{N}$ and of a two-element subgroup of $\sym_{[n]}$ generated by the central symmetry.
\medskip

Let us now consider the task of describing the symmetry group of the permutohedron, but under the {\bf competing convention} that $N:=[n]$, that is, the ground set is identified with the discrete numerical interval $[n]$. What is the difference? The vertices of the permutohedron are now self-bijections $\omega:[n]\to [n]$, that is, the permutations of the set $[n]$. In other words, under this convention, the set of vertices of the permutohedron {\em forms a group}, specifically the symmetric group $\sym_{[n]}$. This was essentially the approach presented by Crisman in \cite{Cri11}; the only technical difference was that he represented the permutations of the set $[n]$ using cycle notation.

Under these circumstances, the elements of the symmetry group of the permutohedron are certain self-transformations of the space $\sym_{[n]}$ and, to describe them, Crisman \cite{Cri11} utilized permutations of $[n]$, that is, individual elements of $\sym_{[n]}$. In other words, he interpreted the elements of $\sym_{[n]}$ as operations acting on the entire set $\sym_{[n]}$, which step then required distinguishing between the so-called ``left" and ``right" actions of this group $\sym_{[n]}$ on itself\,!
Consequently, a multiple ambiguity arose regarding the meaning of a permutation $\omega$ of $[n]$: it can be understood either as a vector in Euclidean space (since $\omega\in\sym_{[n]}\subseteq {\dv R}^{[n]}$) or as one of two possible encodings of a self-transformation of the space $\sym_{[n]}$.
This is certainly difficult for the reader to grasp.

Similar ambiguity issues arose in other publications that adopted the convention $N:=[n]$.
Their authors sometimes attempted to resolve these ambiguities through additional conventions. However,
they usually devoted very little space to explaining their custom supplementary rules,
or even left the task entirely to their readers. For instance, in \cite[\S\,\,2]{MPSSW09},
total orders on $[n]$ were associated with two sorts of vectors: the {\em rank vectors\/} and the {\em descent vectors}. However, the explanation provided there lacks important technical details. The very first sentence
of \cite[\S\,\,2]{MPSSW09} claims that a permutation of $[n]$ is the same as a total order on $[n]$, which is not the case without an additional convention. Indeed, as explained in the beginning of this little essay, there are two different ways to view a permutation of $[n]$ as a total order on $[n]$.
\smallskip

In summary, we have observed that some authors \cite{Zie95,BS96,PRW08,MPSSW09,Cri11,MUWY18} working with the permutohedron adopted a seemingly innocuous  identification convention that its ground set~$N$ coincides with the discrete interval $[n]:=\{1,\ldots,n\}$. They likely viewed this as a harmless convention, identifying the objects of interest with well-known mathematical objects, namely with {\em permutations} of a finite set $[n]$.
As we tried to explain in this essay, such a highly specific convention is unfortunately detrimental, because it leads to misunderstandings on the side of the reader. In our view, perhaps a subjective one, this step amounts to imposing a redundant, irrelevant, and misleading mathematical structure on the set of interest, namely the group structure\,! Indeed, we think that the imposed group structure on the permutohedron gives rise to ambiguities and subsequent identification issues. To avoid these problems in our paper, we came with the conventions from Section~\ref{ssec.notation-convention}.

\end{document}